\documentclass[12pt]{article}
\usepackage[T1]{fontenc}
\usepackage[utf8]{inputenc}
\usepackage{geometry}
\usepackage[dvipsnames]{xcolor}
\usepackage{graphicx}
\usepackage{amsmath,amssymb}
\usepackage{amsmath,mathrsfs}
\usepackage{mathtools}
\usepackage{proof}
\usepackage{amsfonts}
\usepackage{amssymb}
\usepackage{stmaryrd}
\usepackage{pstricks}
\usepackage{setspace}
\usepackage[english]{varioref}
\usepackage{placeins}
\usepackage[english]{babel}
\usepackage{tikz}
\usepackage{makeidx}
\graphicspath{{images/}}
\usepackage[babel=true]{csquotes}
\newcommand{\comment}[1]{}
\newif\ifpdf
\ifpdf \else \fi \textwidth = 6.5 in \textheight = 9 in
\usepackage{bbm}
\usepackage{thmtools}
\usepackage{amsthm}
\usepackage{babel}

\usepackage{mathtools}

\newenvironment{proof}[1][Proof]{\noindent\textbf{#1.} }{$\Box$}

\providecommand{\U}[1]{\protect\rule{.lin}{.lin}}

\newtheorem{theorem}{Theorem}

\newtheorem{corollary}[theorem]{Corollary}

\theoremstyle{definition}

\newtheorem{observation}[theorem]{Observation}

\newtheorem{proposition}[theorem]{Proposition}

\usepackage{caption}
\usepackage{subcaption}

\begin{document}

\title{Quorum coloring of maximum cardinality in linear time for a subclass of perfect trees}

\author{Wissam Boumalha$^{1}$, Asmaa Issad$^{1}$ and Rafik Sahbi$^{2}$\\$^{1}${\small Department of Mathematics, University of Blida 1}\\{\small E-mails: w.boumalha.edu@gmail.com}\\{\small a.issad.edu@gmail.com}\\$^{2}${\small Department of the Preparatory Training}\\{\small Algiers Higher School of Applied Sciences}\\{\small B.P. 474, Martyrs Square, Algiers 16001, Algeria.}\\{\small E-mail: r.sahbi@g.essa-alger.edu.dz}}
\date{}
\maketitle

\begin{abstract}
A partition $\pi=\{V_{1},V_{2},...,V_{k}\}$ of the vertex set $V$ of a graph $G$ into $k$ color classes $V_{i}$, with $1\leq i\leq k$ is called a quorum coloring of $G$ if for every vertex $v\in V$, at least half of the vertices in the closed neighborhood $N[v]$ of $v$ have the same color as $v$. The maximum cardinality of a quorum coloring of $G$ is called the quorum coloring number of $G$ and is denoted by $\psi_{q}(G)$. A quorum coloring of order $\psi_{q}(G)$ is a $\psi_{q}$-coloring. The determination of the quorum coloring number or design a linear-time algorithm computing it in a perfect $N$-ary tree has been posed recently as an open problem by Sahbi [R. Sahbi, New sharp lower bound for the quorum coloring number of trees, Inf. Process. Lett. 178, (2022) 106297]. In this paper, we answer this problem by designing a linear-time algorithm for finding both a $\psi_{q}$-coloring and the quorum coloring number for every perfect tree whose the vertices at the same depth have the same degree.
\end{abstract}

\section{Introduction}

\subsection{Definitions and notations}

All the graphs of this paper are simple, that is, they are finite, undirected and have neither loops nor multiple edges.

Let $G=(V,E)$ be a graph. The {\it order} of $G$ is denoted by $n$.The {\it induced subgraph} of $G$ by a subset $S$ of $V$ is denoted by $G[S]$. For every vertex $v\in V,$ the {\it open neighborhood} $N_{G}(v)$ is the set $\{u\in V(G):uv\in E(G)\}$ and the {\it closed neighborhood} of $v$ is the set $N_{G}[v]=N_{G}(v)\cup\{v\}$. The {\it degree} of a vertex $v\in V(G)$ in $G$ is $|N_{G}(v)|$. A vertex of degree zero in $G$ is an {\it isolated} vertex of $G$ and a vertex of degree one in $G$ is a {\it leaf} or a {\it pendent vertex} of $G$. The set of leaves of $G$ is denoted by $L(G)$, or simply $L$ when $G$ is unambiguous. The maximum and minimum vertex degrees in $G$ are respectively denoted by $\Delta(G)$ and $\delta(G)$. More generally for a vertex $v\in V$ and a subset $S\subseteq V$, the {\it open} and {\it closed neighborhoods} of $v$ in $S$ are respectively defined by the sets $N_{S}(v)=\{u\in S:uv\in E(G)\}$ and $N_{S}[v]=N_{S}(v)\cup\{v\}$, and the degree of $v$ in $S$ is $d_{S}(v)=|N_{S}(v)|$. In particular, one can see that for $S=V$ we have $N_{S}(v)=N_{G}(v)$, $N_{G}[v]=N_{G}[v]$ and $d_{S}(v)=d_{G}(v)$. A {\it tree} is a connected graph having no cycle. In this paper, we denote by {\it binary tree} every tree $T$ with $\Delta(T)\leq 3$. The {\it distance} $d_{G}(u,v)$ between two vertices $u$ and $v$ in a connected graph $G$ is the length of the shortest $u-v$ path in $G$, and the {\it diameter} of $G$ is the longest distance $\displaystyle diam(G)=\max_{u,v\in V}d_{G}(u,v)$ between two vertices in $G$. A {\it matching} in a graph $G=(V,E)$ is a set $M\subseteq E$ whose edges are pairwise non adjacent. The {\it matching number} $\mu(G)$ equals the maximum cardinality of a matching in $G.$

A {\it rooted} tree $T=(V,E,r)$ is a tree $T=(V,E)$ with a distinguished vertex $r$ called the {\it root} of $T$. A vertex $v$ is the {\it parent} of a vertex $w$ in a rooted tree $T=(V,E,r)$ if $vw\in E$ and $d_{T}(w,r)=d_{T}(v,r)+1,$ in which case $w$ is said to be a {\it child} of $v$. Two vertices having the same parent in a rooted tree $T=(V,E,r)$ are called {\it siblings}. A vertex $x$ is a {\it descendant} of a vertex $v$ in a rooted tree $T=(V,E,r)$ if $d_{T}(v,x)=d_{T}(x,r)-d_{T}(v,r)$ (in particular, we have $d_{T}(x,r)>d_{T}(v,r)$); in this case, $v$ is said to be an {\it ancestor} of $x$. 

A {\it perfect} tree is a rooted tree $T=(V,E,r)$ whose leaves are at the same distance $h$ from the root $r$, distance called the {\it height} of $T$ and defined by $\displaystyle h=\max_{v\in V}d_{T}(r,v)$. Also, we denote by $D_{i}$ the set $\{v\in V~|~d_{T}(r,v)=i\}$ of vertices of a perfect tree $T=(V,E,r)$ that are at distance $i$ from $r$, distance called the {\it depth} of $D_{i}$'s vertices, for every $i\in\{0,1,\ldots,h\}$. We also set $\displaystyle D_{i}=\{v_{i,j}\}_{1\leq j\leq\ell_{i}}$ for every $i\in\{0,1,\ldots,h\}$ (in particular, $r=v_{0,1}$). 

A {\it perfect} $N_{i}${\it -ary} tree {\it per level} is a perfect tree $T=(V,E,r)$ whose all the $D_{i}$'s vertices have the same degree $N_{i}+1$, that is, all the vertices that are at the same depth $i$ have the same number of children $N_{i}$. It can be seen from the above that every siblings of a perfect $N_{i}$-ary tree by level have the same degree. In particular, if $N_{i}=N$ for some integer $N\geq 2$ and every $\displaystyle i\in\left\{0,\ldots,h-1\right\}$ and $j\in\{1,\ldots,\ell_{i}\}$, then $T$ is called a {\it perfect} $N${\it -ary tree} (for $N=2$, $T$ is called {\it perfect binary} tree). Note from the definitions that every perfect $N$-ary tree is a perfect $N_{i}$-ary tree per level with $N_{i}=N$ for every $i\in\{0,1,\ldots,h-1\}$, but that obviously, the converse is not true.

A partition $\pi=\{V_{1},V_{2},...,V_{k}\}$ of the vertex set $V$ of a graph $G$ into $k$ color classes $V_{i}$, with $i\in\{1,...,k\}$ is called a \textit{quorum coloring} of $G$ if for every vertex $v\in V,$ at least half of the vertices in the closed neighborhood $N_{G}[v]$ have the same color as $v$, which means formally that for every $i\in\{1,\ldots,k\}$ and every $v\in V_{i}$, we have $|N_{G}[v]\cap V_{i}|\geq\frac{|N_{G}[v]|}{2}.$ The color classes $V_{i}$ are called {\it quorum classes} and for a fixed $i\in\{1,\ldots,k\}$, each vertex of $V_{i}$ is called a {\it quorum vertex}. The maximum cardinality of a quorum coloring of $G$ is called the \textit{quorum coloring number} of $G$ and is denoted by $\psi_{q}(G)$. A quorum coloring of $G$ of cardinality $\psi_{q}(G)$ is called a $\psi_{q}$-\textit{coloring} of $G$. For a quorum coloring $\pi$ on a graph $G$ and any vertex $v\in V(G)$, we denote the class of $v$ with respect to $\pi$ by $\pi(v)$.

\subsection{Previous results}

Quorum colorings is another name of partitions into defensive alliances in graphs (cf. \cite{EroAll}) that admit applications to data clustering (see \cite{OuazSurv,SahSol,ShaThes}), and the quorum classes of a graph with respect to a given quorum coloring are nothing but defensive alliances. For further reading on the well-studied topic of defensive alliances in graphs, the reader can consult the references \cite{FERsur,Fricke,HayHedChapter,HayAll,KriAll,KRISintro,OuazSurv,ShaThes,YeroSurv,
YeroPart}.

Quorum colorings in graphs were introduced in 2013 by Hedetniemi et al. in \cite{HedQuor} where the authors studied their basic properties as well as the quorum coloring number of some usual graphs among which the hypercubes, then concluded their article by raising and listing twelve open problems. In 2018, Sahbi and Chellali \cite{SahQuor} answered three of these open problems and showed in particular that the decision problem associated with $\psi_{q}(G)$ is {\it NP}-complete for a general graph $G$. Furthermore, Sahbi pursued the study of the open questions raised by the authors \cite{HedQuor} throughout \cite{SahNew,SahComp,SahSol} and in \cite{SubQuor} with Belkina and Bennadji. In particular, he brought a partial answer to the following open question posed by Hedetniemi et al.

\begin{description}

\item[1.] Can you design a linear-time algorithm for computing the value of $\psi_{q}(T)$ for any tree $T$ ?

\end{description}

In fact, the author \cite{SahNew} first established a lower bound on the quorum coloring number of nontrivial trees involving the matching number of the subtree induced by the non pendent vertices, the order of the tree and the vertex degrees and showed that this bound is computable in linear time. 

\begin{theorem}\label{The0}\cite{SahNew} Let $T=(V,E)$ be a nontrivial tree. Then, \[\psi_{q}(T)\geq\mu\left(T[V\setminus L] \right)+n-\displaystyle\sum_{v\in\left(V\setminus L\right)}\left\lfloor\dfrac{d_{T}(v)}{2}\right\rfloor.\] This bound can be computed in linear time.
\end{theorem}

Then, he showed that the bound of Theorem \ref{The0} is attained by all binary trees. 


\begin{corollary}\label{Cor0}\cite{SahNew} Let $T=(V,E)$ be a non trivial binary tree. Then, \[\psi_{q}(T)=\mu(T[V\setminus L])+|L|.\]
This value can be computed in linear time.
\end{corollary}

Finally, Sahbi \cite{SahNew} posed the following open problem whose the second part is none other than {\it Question 1} restricted to perfect $N$-ary trees.

\begin{description}

\item[2.] Determine the exact value of the quorum coloring number of perfect $N$-ary trees or design a linear-time algorithm computing it.

\end{description}

In this paper, we first recall some fundamental results on quorum colorings of graphs in Section 2, then we answer the second part of Problem 2 in Section 3 by designing a linear-time algorithm both for finding a $\psi_{q}$-coloring and computing the quorum coloring number of any perfect $N_{i}$-ary tree per level.

\section{Preliminary results}

This section is devoted to the statement of some fundamental results on quorum colorings that have been previously established in the literature, some of which will be used in this work. We start by a linear relationship between the quorum coloring number of a disconnected graph and those of its components.

\begin{proposition}\label{Pro1}\cite{EroAll} Let $G$ be a disconnected graph whose components are $G_{1},G_{2},\ldots,G_{r}$ ($r\geq1$). Then \[\displaystyle\psi_{q}(G)=\sum_{1\leq i\leq r}\psi_{q}(G_{i}).\]
\end{proposition}

The next result provides sharp lower and upper bounds of the quorum coloring number of a disconnected graph of order at least three. However, one can easily see that this result remains true for any graph and any order.

\begin{proposition}\label{Pro2}\cite{EroAll} Let $G$ be a disconnected graph of order $n\geq3$. Then \[1\leq\psi_{q}(G)\leq n.\]\end{proposition}

The sharpness of the bounds of Proposition \ref{Pro2} can be seen thanks to the following two propositions.

\begin{proposition}\label{Pro3}\cite{HedQuor} For the complete graph $K_{n}$ of odd order, $\psi_{q}(K_{n})=1$, while for any complete graph $K_{n}$ of even order, $\psi_{q}(K_{n})=2$.\end{proposition}

\begin{proposition}\label{Pro4}\cite{HedQuor} Let $G$ be a graph of order $n$. Then $\psi_{q}(G)=n$ if and only if $\Delta(G)\leq 1$, that is $G$ consists of isolated vertices and disjoint copies of the complete graph $K_{2}$ of order $2$.\end{proposition}

The next result was proved by Hedetniemi et al. It says that in a $\psi_{q}$-coloring, a quorum class always induces a connected subgraph.

\begin{proposition}\label{Pro5}\cite{HedQuor} Let $G$ be a graph, and let $\pi=\{V_{1},V_{2},\ldots,V_{k}\}$ be any $\psi_{q}$-coloring of $G$. Then, for every $i$, $1\leq i\leq k$, the induced subgraph $G[V_{i}]$ is connected.
\end{proposition}

In \cite{SahNew}, Sahbi proved the following observation which gives four equivalent assertions to the quorum vertex property and the second of which will be used in section 3.

\begin{observation}\label{Obs6}\cite{SahNew} Let $G=(V,E)$ be a graph, $\pi=\{V_{1},V_{2},\ldots,V_{k}\}$ a quorum coloring of $G$ and $i\in\{1,2,\ldots,k\}$. Then the following assertions are equivalent:
\begin{itemize}
\item[1.] $V_{i}$ is a quorum class.
\item[2.] For every vertex $v\in V_{i},~|N_{G}[v]\cap V_{i}|\geq\left\lceil\frac{\left|N_{G}[v]\right|}{2}\right \rceil$.
\item[3.] For every vertex $v\in V_{i},~|N_{G}[v]\cap V_{i}|\geq|N_{G}[v]\cap(V\setminus V_{i})|$.
\item[4.] For every vertex $v\in V_{i},~\displaystyle d_{V_{i}}(v)+1\geq d_{V\setminus V_{i}}(v)$.
\item[5.] For every vertex $v\in V_{i},~\displaystyle d_{V_{i}}(v)\geq\left\lfloor\dfrac{d_{G}(v)}{2}\right \rfloor$.
\end{itemize}
\end{observation}

As consequences of Observation \ref{Obs6}, the author \cite{SahNew} deduced the three following corollaries.

\begin{corollary}\label{Cor7}\cite{SahNew} Let $G=(V,E)$, $\pi=\{V_{1},V_{2},\ldots,V_{k}\}$ a quorum coloring of $G$ and $i\in\{1,2,\ldots,k\}$. Then for every vertex $v\in V_{i}$, $|V_{i}|\geq\left\lfloor\dfrac{d_{G}(v)}{2}\right\rfloor+1$.
\end{corollary}

\begin{corollary}\label{Cor8}\cite{SahNew} Let $G=(V,E)$ be a graph, $\pi=\{V_{1},V_{2},\ldots,V_{k}\}$ a $\psi_{q}$-coloring of $G$, $i\in\{1,2,\ldots,k\}$ a positive integer and $v\in V_{i}$ a vertex such that $\displaystyle d_{G}(v)=\max_{u\in V_{i}}d_{G}(u)$. Then $|V_{i}|=1$ if and only if $d_{G}(v)\leq 1$.
\end{corollary}

\begin{corollary}\label{Cor9}\cite{SahNew} Let $G=(V,E)$, $\pi=\{V_{1},V_{2},\ldots,V_{k}\}$ a $\psi_{q}$-coloring of $G$ and $i\in\{1,2,\ldots,k\}$ a positive integer. Then $|V_{i}|\geq2$ if and only if $V_{i}$ contains a vertex of degree at least $2$.
\end{corollary}

Corollary \ref{Cor7} provides a lower bound of a quorum class in term of maximum degree of its vertices, while Corollary \ref{Cor8} states that the unique vertex of a singleton quorum class is necessarily pendent or isolated. Corollary \ref{Cor9}, obtained by negating Corollary \ref{Cor8}, shows that every vertex of degree at least two is contained in a quorum class of order at least two in any $\psi_{q}$-coloring, but this result can easily be extended to a quorum coloring that is not necessarily a $\psi_{q}$-coloring.

In the next section, we prove our result announced in Section 1.

\section{Answer to Problem 2}

In this section, we answer Problem 2 stated in Section 1. To do this, we need the following definitions.

A quorum coloring $\pi$ of a graph $G$ is said to be {\it cost-effective} if every vertex of $V\setminus L$ satisfies the quorum vertex property with equality, that is, if for every vertex $v\in V$ we have $|N_{G}[v]\cap\pi(v)|=\left\lceil\frac{|N_{G}[v]|}{2}\right\rceil$. Given two quorum colorings $\pi_{1}$ and $\pi_{2}$ of $G$, we say that $\pi_{1}$ is {\it better} than $\pi_{2}$ if $\pi_{1}$ is cost-effective and $|\pi_{1}|\geq |\pi_{2}|$.

Before proving our main result, we first prove that we can modify any quorum coloring of an arbitrarily tree to obtain a better one in the sense of the above definition.

\begin{theorem}\label{QC->CEF} Let $T=(V,E,r)$ be a tree and $\pi_{0}$ a quorum coloring of $T$. Then, there exists a quorum coloring of $T$ that is better than $\pi_{0}$.
\end{theorem}

\begin{proof} Let $\pi_{0}$ be a quorum coloring of $T$. We run the following algorithm that we denote by Algorithm 1.

\begin{itemize}

\item[$\bullet$] For $i=0$ to $h-1$ and $j=1$ to $\ell_{i}$, execute the following two steps.

\begin{description}

\item[1.] For $|N_{T}[v_{i,j}]\cap\pi_{i}(v_{i,j})|>\left\lceil\frac{|N_{T}[v_{i,j}]|}{2}\right\rceil$, chose arbitrarily a subset $S_{i,j}$ of $\pi_{i}(v_{i,j})\cap D_{i+1}$ so that $|N_{T}[v_{i,j}]\cap\pi_{i}(v_{i,j})|-|S_{i,j}|=\left\lceil\frac{|N_{T}[v_{i,j}]|}{2}\right\rceil$, color the $S_{i,j}$'s vertices with $|S_{i,j}|$ new colors and for every vertex $v\in S_{i,j}$, assign the new $v$'s color to all the descendants of $v$ that were in $\pi_{i}(v)$; we denote the obtained coloring by $\pi_{i}'$. Then, go to step 2.

\item[2.] For every vertex $v\in S_{i,j}$ such that $|N_{T}[v]\cap\pi_{i}'(v)|<\left\lceil\frac{|N_{T}[v]|}{2}\right\rceil$ (note in this case that we have necessarily $|N_{T}[v]\cap\pi_{i}'(v)|=\left\lceil\frac{|N_{T}[v]|}{2}\right\rceil-1$ since by passing from $\pi_{i}$ to $\pi_{i}'$, the unique neighbor of $v$ that has no longer the $v$'s color is $v_{i,j}$ implying that $|\pi_{i}(v)|-|\pi_{i}'(v)|=1$), assign the $v$'s color to an arbitrarily chosen vertex $w$ from $\left(N_{T}(v)\cap D_{i+2}\right)\setminus\pi_{i}'(v)$ and to all the descendants of $w$ that belonged to $\pi_{i}'(v)$; therefore, it can be seen without difficulty that the resulted coloring is a quorum coloring that we denote by $\pi_{i+1}$.

\end{description}

\end{itemize}

Algorithm 1 must terminate since $T$ is of finite order. When it is finished, we clearly obtain a cost-effective quorum coloring $\pi_{h}$ according to steps 1 and 2. Moreover, each time Algorithm 1 is run for some $i\in\{0,\ldots,h-1\}$ and $j\in\{1,\ldots,\ell_{i}\}$, the number of colors increases by at least one at the execution of step 1, while it decreases by at most one at the execution of step 2 and consequently, the number of colors does not decrease. Hence we deduce that $\pi_{h}$ is better than $\pi_{0}$. \end{proof}

Starting from a quorum coloring $\pi_{0}$ of any tree $T$, the proof of Theorem \ref{QC->CEF} shows that one can obtain a quorum coloring of $T$ better than $\pi_{0}$ that is cost-effective. As consequence, if $\pi_{0}$ is a $\psi_{q}$-coloring then Algorithm 1 concludes with a cost effective $\psi_{q}$-coloring of $T$ as output.

\begin{corollary}\label{corQC->CEF}
Every tree has a cost-effective $\psi_{q}$-coloring.\end{corollary}

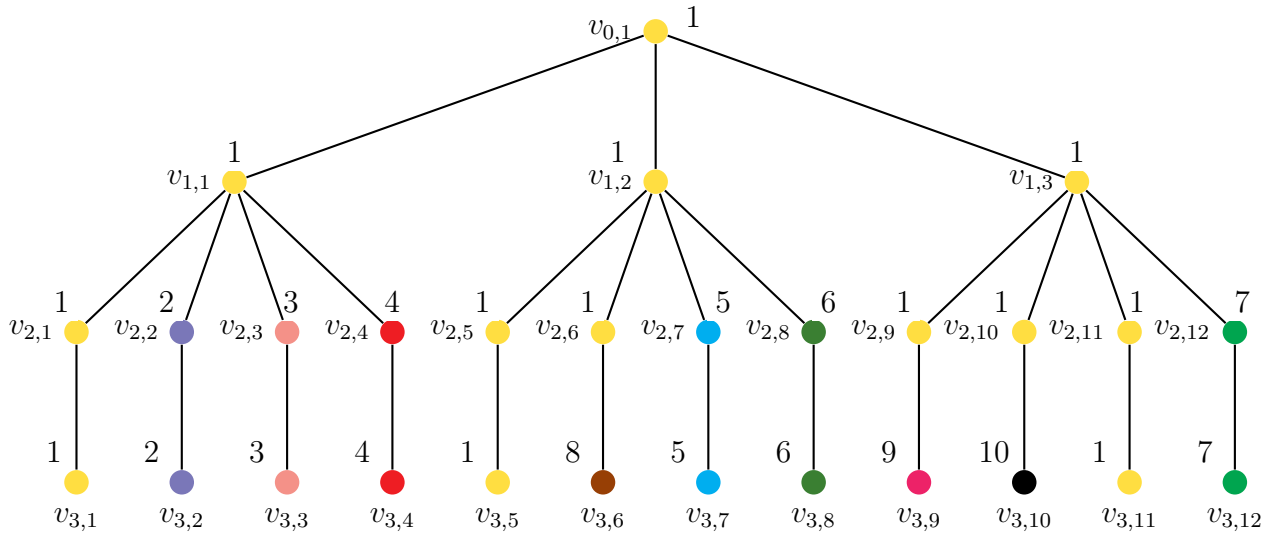
\begin{figure}[ht!]
\begin{center}
\begin{tikzpicture}[inner sep=1.15mm]
\tikzstyle{a}=[circle,fill=Periwinkle!]
\tikzstyle{A}=[circle,fill=Salmon!]
\tikzstyle{L}=[circle,fill=YellowOrange!]
\tikzstyle{b}=[circle,fill=Red!100]
\tikzstyle{o}=[circle,draw,fill=White!]
\tikzstyle{M}=[circle,fill=Dandelion!]
\tikzstyle{O}=[circle,fill=Goldenrod!]
\tikzstyle{S}=[circle,fill=Orange!100]
\tikzstyle{m}=[circle,fill=SeaGreen!]
\tikzstyle{d}=[circle,fill=OliveGreen!]
\tikzstyle{D}=[circle,fill=WildStrawberry!]
\tikzstyle{Q}=[circle,fill=LimeGreen!]
\tikzstyle{N}=[circle,fill=Thistle!]
\tikzstyle{p}=[circle,fill=Fuchsia!]
\tikzstyle{P}=[circle,fill=RoyalPurple!]
\tikzstyle{R}=[circle,fill=Magenta!]
\tikzstyle{I}=[circle,fill=Yellow!100]
\tikzstyle{i}=[circle,fill=GreenYellow!]
\tikzstyle{K}=[circle,fill=CarnationPink!]
\tikzstyle{w}=[circle,fill=Bittersweet!]
\tikzstyle{E}=[circle,fill=black!100]
\tikzstyle{n}=[rectangle,fill=black!0]
\tikzstyle{h}=[circle,fill=blue!100]
\tikzstyle{s}=[circle,fill=TealBlue!]
\tikzstyle{H}=[circle,fill=RoyalBlue!]
\tikzstyle{C}=[circle,fill=Cyan!]
\tikzstyle{f}=[circle,fill=Green!]
\tikzstyle{g}=[circle,fill=SpringGreen!]
\tikzstyle{G}=[circle,fill=Tan!]
\tikzstyle{r}=[circle,fill=Gray!]
\tikzstyle{c}=[circle,fill=CadetBlue!]
\tikzstyle{j}=[circle,fill=pink!]
\tikzstyle{k}=[circle,fill=Peach!]
\tikzstyle{q}=[circle,fill=JungleGreen!]
\tikzstyle{t}=[circle,fill=Apricot!]
\tikzstyle{l}=[circle,fill=BrickRed!]
\tikzstyle{B}=[circle,fill=RawSienna!]
\tikzstyle{J}=[circle,fill=Turquoise!40]
\tikzstyle{F}=[circle,fill=Black!20]
\tikzstyle{T}=[circle,fill=Black!40]
\tikzstyle{u}=[circle,fill=BrickRed!40]
\tikzstyle{U}=[circle,fill=black!50]
\tikzstyle{1}=[circle,fill=brown!100]
\tikzstyle{2}=[circle,fill=purple!100]
\tikzstyle{3}=[circle,fill=green!50]
\tikzstyle{4}=[circle,fill=orange!]
\tikzstyle{e}=[-,thick]
\node [O](v1)at (0,0){};\node [n](v1111)at (0,-0.5){$v_{3,1}$};
\node [a](v2)at (1.4,0){};\node [n](v2222)at (1.4,-0.5){$v_{3,2}$};
\node [A](v3)at (2.8,0){};\node [n](v3333)at (2.8,-0.5){$v_{3,3}$};
\node [b](v4)at (4.2,0){};\node [n](v4444)at (4.2,-0.5){$v_{3,4}$};
\node [O](v5)at (5.6,0){};\node [n](v5555)at (5.6,-0.5){$v_{3,5}$};
\node [B](v6)at (7,0){};\node [n](v6666)at (7,-0.5){$v_{3,6}$};
\node [C](v7)at (8.4,0){};\node [n](v7777)at (8.4,-0.5){$v_{3,7}$};
\node [d](v8)at (9.8,0){};\node [n](v8888)at (9.8,-0.5){$v_{3,8}$};
\node [D](v9)at (11.2,0){};\node [n](v9999)at (11.2,-0.5){$v_{3,9}$};
\node [E](v10)at (12.6,0){};\node [n](v101010)at (12.6,-0.5){$v_{3,10}$};
\node [O](v11)at (14,0){};\node [n](v111111)at (14,-0.5){$v_{3,11}$};
\node [f](v12)at (15.4,0){};\node [n](v121212)at (15.4,-0.5){$v_{3,12}$};
\node [O](v13)at (0,2){};\node [n](v131313)at (-0.6,2){$v_{2,1}$};
\node [a](v14)at (1.4,2){};\node [n](v1414)at (0.8,2){$v_{2,2}$};
\node [A](v15)at (2.8,2){};\node [n](v1515)at (2.2,2){$v_{2,3}$};
\node [b](v16)at (4.2,2){};\node [n](v1616)at (3.6,2){$v_{2,4}$};
\node [O](v17)at (5.6,2){};\node [n](v1717)at (5,2){$v_{2,5}$};
\node [O](v18)at (7,2){};\node [n](v1818)at (6.4,2){$v_{2,6}$};
\node [C](v19)at (8.4,2){};\node [n](v1919)at (7.8,2){$v_{2,7}$};
\node [d](v20)at (9.8,2){};\node [n](v2020)at (9.2,2){$v_{2,8}$};
\node [O](v21)at (11.2,2){};\node [n](v2121)at (10.6,2){$v_{2,9}$};
\node [O](v22)at (12.6,2){};\node [n](v22222)at (11.9,2){$v_{2,10}$};
\node [O](v23)at (14,2){};\node [n](v2323)at (13.3,2){$v_{2,11}$};
\node [f](v24)at (15.4,2){};\node [n](v2424)at (14.7,2){$v_{2,12}$};
\node [O](v25)at (2.1,4){};\node [n](v2525)at (1.5,4){$v_{1,1}$};
\node [O](v26)at (7.7,4){};\node [n](v2626)at (7.1,4){$v_{1,2}$};
\node [O](v27)at (13.3,4){};\node [n](v2727)at (12.7,4){$v_{1,3}$};
\node [O](v28)at (7.7,6){};\node [n](v2828)at (7.1,6){$v_{0,1}$};
\node [n](v11111)at (-0.3,0.4){$1$};
\node [n](v22222)at (1,0.4){$2$};
\node [n](v33333)at (2.4,0.4){$3$};
\node [n](v44444)at (3.8,0.4){$4$};
\node [n](v55555)at (5.2,0.4){$1$};
\node [n](v66666)at (6.6,0.4){$8$};
\node [n](v77777)at (8,0.4){$5$};
\node [n](v88888)at (9.4,0.4){$6$};
\node [n](v99999)at (10.8,0.4){$9$};
\node [n](v10101010)at (12.2,0.4){$10$};
\node [n](v11111111)at (13.6,0.4){$1$};
\node [n](v12121212)at (15,0.4){$7$};
\node [n](v13131313)at (-0.2,2.4){$1$};
\node [n](v141414)at (1.2,2.4){$2$};
\node [n](v151515)at (2.85,2.4){$3$};
\node [n](v161616)at (4.2,2.4){$4$};
\node [n](v171717)at (5.4,2.4){$1$};
\node [n](v181818)at (6.8,2.4){$1$};
\node [n](v191919)at (8.6,2.4){$5$};
\node [n](v202020)at (10,2.4){$6$};
\node [n](v212121)at (11,2.4){$1$};
\node [n](v2222222)at (12.3,2.4){$1$};
\node [n](v232323)at (14.1,2.4){$1$};
\node [n](v242424)at (15.5,2.4){$7$};
\node [n](v252525)at (2.1,4.4){$1$};
\node [n](v262626)at (7.2,4.4){$1$};
\node [n](v272728)at (13.3,4.4){$1$};
\node [n](v282828)at (8.2,6.2){$1$};
\draw[e](v1)--(v13);\draw[e](v2)--(v14);\draw[e](v3)--(v15);
\draw[e](v4)--(v16);\draw[e](v5)--(v17);\draw[e](v6)--(v18);\draw[e](v7)--(v19);\draw[e](v8)--(v20);\draw[e](v9)--(v21);\draw[e](v10)--(v22);\draw[e](v11)--(v23);\draw[e](v12)--(v24);\draw[e](v13)--(v25);\draw[e](v14)--(v25);\draw[e](v15)--(v25);\draw[e](v16)--(v25);\draw[e](v17)--(v26);\draw[e](v18)--(v26);\draw[e](v19)--(v26);\draw[e](v20)--(v26);\draw[e](v21)--(v27);\draw[e](v22)--(v27);\draw[e](v23)--(v27);\draw[e](v24)--(v27);\draw[e](v28)--(v27);\draw[e](v28)--(v26);\draw[e](v25)--(v28);
\end{tikzpicture}
\caption{An arbitrary initial quorum coloring $\pi_{0}$ with $|\pi_{0}|=10$}
\label{1}
\end{center}
\end{figure}

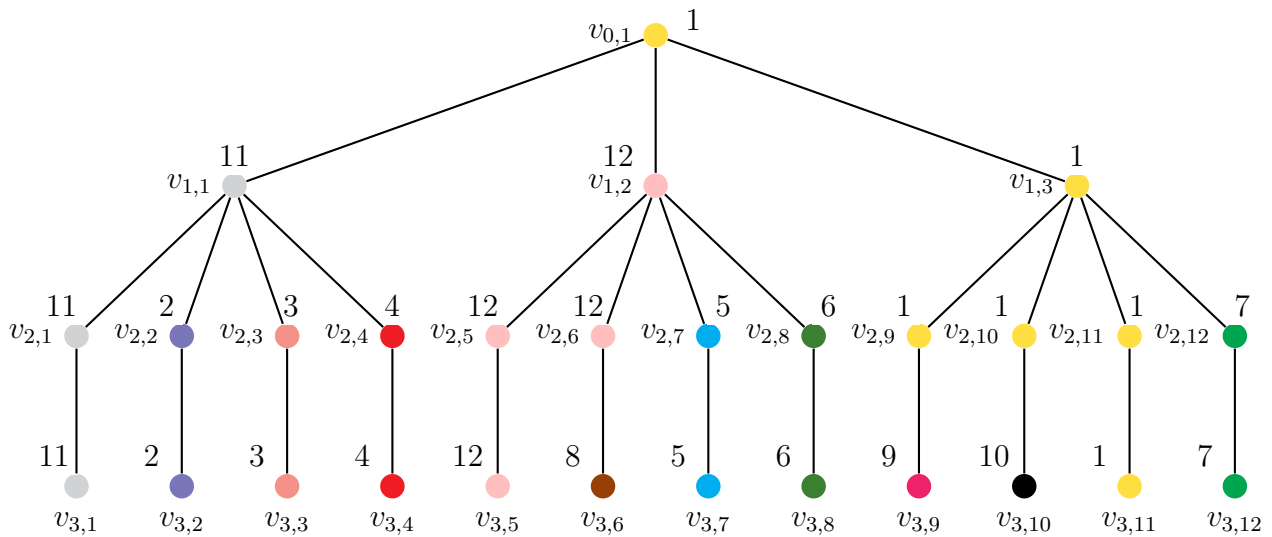
\begin{figure}[ht!]
\begin{center}
\begin{tikzpicture}[inner sep=1.15mm]
\tikzstyle{a}=[circle,fill=Periwinkle!]
\tikzstyle{A}=[circle,fill=Salmon!]
\tikzstyle{L}=[circle,fill=YellowOrange!]
\tikzstyle{b}=[circle,fill=Red!100]
\tikzstyle{o}=[circle,draw,fill=White!]
\tikzstyle{M}=[circle,fill=Dandelion!]
\tikzstyle{O}=[circle,fill=Goldenrod!]
\tikzstyle{S}=[circle,fill=Orange!100]
\tikzstyle{m}=[circle,fill=SeaGreen!]
\tikzstyle{d}=[circle,fill=OliveGreen!]
\tikzstyle{D}=[circle,fill=WildStrawberry!]
\tikzstyle{Q}=[circle,fill=LimeGreen!]
\tikzstyle{N}=[circle,fill=Thistle!]
\tikzstyle{p}=[circle,fill=Fuchsia!]
\tikzstyle{P}=[circle,fill=RoyalPurple!]
\tikzstyle{R}=[circle,fill=Magenta!]
\tikzstyle{I}=[circle,fill=Yellow!100]
\tikzstyle{i}=[circle,fill=GreenYellow!]
\tikzstyle{K}=[circle,fill=CarnationPink!]
\tikzstyle{w}=[circle,fill=Bittersweet!]
\tikzstyle{E}=[circle,fill=black!100]
\tikzstyle{n}=[rectangle,fill=black!0]
\tikzstyle{h}=[circle,fill=blue!100]
\tikzstyle{s}=[circle,fill=TealBlue!]
\tikzstyle{H}=[circle,fill=RoyalBlue!]
\tikzstyle{C}=[circle,fill=Cyan!]
\tikzstyle{f}=[circle,fill=Green!]
\tikzstyle{g}=[circle,fill=SpringGreen!]
\tikzstyle{G}=[circle,fill=Tan!]
\tikzstyle{r}=[circle,fill=Gray!]
\tikzstyle{c}=[circle,fill=CadetBlue!]
\tikzstyle{j}=[circle,fill=pink!]
\tikzstyle{k}=[circle,fill=Peach!]
\tikzstyle{q}=[circle,fill=JungleGreen!]
\tikzstyle{t}=[circle,fill=Apricot!]
\tikzstyle{l}=[circle,fill=BrickRed!]
\tikzstyle{B}=[circle,fill=RawSienna!]
\tikzstyle{J}=[circle,fill=Turquoise!40]
\tikzstyle{F}=[circle,fill=Black!20]
\tikzstyle{T}=[circle,fill=Black!40]
\tikzstyle{u}=[circle,fill=BrickRed!40]
\tikzstyle{U}=[circle,fill=black!50]
\tikzstyle{1}=[circle,fill=brown!100]
\tikzstyle{2}=[circle,fill=purple!100]
\tikzstyle{3}=[circle,fill=green!50]
\tikzstyle{4}=[circle,fill=orange!]
\tikzstyle{e}=[-,thick]
\node [F](v1)at (0,0){};\node [n](v1111)at (0,-0.5){$v_{3,1}$};
\node [a](v2)at (1.4,0){};\node [n](v2222)at (1.4,-0.5){$v_{3,2}$};
\node [A](v3)at (2.8,0){};\node [n](v3333)at (2.8,-0.5){$v_{3,3}$};
\node [b](v4)at (4.2,0){};\node [n](v4444)at (4.2,-0.5){$v_{3,4}$};
\node [j](v5)at (5.6,0){};\node [n](v5555)at (5.6,-0.5){$v_{3,5}$};
\node [B](v6)at (7,0){};\node [n](v6666)at (7,-0.5){$v_{3,6}$};
\node [C](v7)at (8.4,0){};\node [n](v7777)at (8.4,-0.5){$v_{3,7}$};
\node [d](v8)at (9.8,0){};\node [n](v8888)at (9.8,-0.5){$v_{3,8}$};
\node [D](v9)at (11.2,0){};\node [n](v9999)at (11.2,-0.5){$v_{3,9}$};
\node [E](v10)at (12.6,0){};\node [n](v101010)at (12.6,-0.5){$v_{3,10}$};
\node [O](v11)at (14,0){};\node [n](v111111)at (14,-0.5){$v_{3,11}$};
\node [f](v12)at (15.4,0){};\node [n](v121212)at (15.4,-0.5){$v_{3,12}$};
\node [F](v13)at (0,2){};\node [n](v131313)at (-0.6,2){$v_{2,1}$};
\node [a](v14)at (1.4,2){};\node [n](v1414)at (0.8,2){$v_{2,2}$};
\node [A](v15)at (2.8,2){};\node [n](v1515)at (2.2,2){$v_{2,3}$};
\node [b](v16)at (4.2,2){};\node [n](v1616)at (3.6,2){$v_{2,4}$};
\node [j](v17)at (5.6,2){};\node [n](v1717)at (5,2){$v_{2,5}$};
\node [j](v18)at (7,2){};\node [n](v1818)at (6.4,2){$v_{2,6}$};
\node [C](v19)at (8.4,2){};\node [n](v1919)at (7.8,2){$v_{2,7}$};
\node [d](v20)at (9.8,2){};\node [n](v2020)at (9.2,2){$v_{2,8}$};
\node [O](v21)at (11.2,2){};\node [n](v2121)at (10.6,2){$v_{2,9}$};
\node [O](v22)at (12.6,2){};\node [n](v22222)at (11.9,2){$v_{2,10}$};
\node [O](v23)at (14,2){};\node [n](v2323)at (13.3,2){$v_{2,11}$};
\node [f](v24)at (15.4,2){};\node [n](v2424)at (14.7,2){$v_{2,12}$};
\node [F](v25)at (2.1,4){};\node [n](v2525)at (1.5,4){$v_{1,1}$};
\node [j](v26)at (7.7,4){};\node [n](v2626)at (7.1,4){$v_{1,2}$};
\node [O](v27)at (13.3,4){};\node [n](v2727)at (12.7,4){$v_{1,3}$};
\node [O](v28)at (7.7,6){};\node [n](v2828)at (7.1,6){$v_{0,1}$};
\node [n](v11111)at (-0.3,0.4){$11$};
\node [n](v22222)at (1,0.4){$2$};
\node [n](v33333)at (2.4,0.4){$3$};
\node [n](v44444)at (3.8,0.4){$4$};
\node [n](v55555)at (5.2,0.4){$12$};
\node [n](v66666)at (6.6,0.4){$8$};
\node [n](v77777)at (8,0.4){$5$};
\node [n](v88888)at (9.4,0.4){$6$};
\node [n](v99999)at (10.8,0.4){$9$};
\node [n](v10101010)at (12.2,0.4){$10$};
\node [n](v11111111)at (13.6,0.4){$1$};
\node [n](v12121212)at (15,0.4){$7$};
\node [n](v13131313)at (-0.2,2.4){$11$};
\node [n](v141414)at (1.2,2.4){$2$};
\node [n](v151515)at (2.85,2.4){$3$};
\node [n](v161616)at (4.2,2.4){$4$};
\node [n](v171717)at (5.4,2.4){$12$};
\node [n](v181818)at (6.8,2.4){$12$};
\node [n](v191919)at (8.6,2.4){$5$};
\node [n](v202020)at (10,2.4){$6$};
\node [n](v212121)at (11,2.4){$1$};
\node [n](v2222222)at (12.3,2.4){$1$};
\node [n](v232323)at (14.1,2.4){$1$};
\node [n](v242424)at (15.5,2.4){$7$};
\node [n](v252525)at (2.1,4.4){$11$};
\node [n](v262626)at (7.2,4.4){$12$};
\node [n](v272728)at (13.3,4.4){$1$};
\node [n](v282828)at (8.2,6.2){$1$};
\draw[e](v1)--(v13);\draw[e](v2)--(v14);\draw[e](v3)--(v15);
\draw[e](v4)--(v16);\draw[e](v5)--(v17);\draw[e](v6)--(v18);\draw[e](v7)--(v19);\draw[e](v8)--(v20);\draw[e](v9)--(v21);\draw[e](v10)--(v22);\draw[e](v11)--(v23);\draw[e](v12)--(v24);\draw[e](v13)--(v25);\draw[e](v14)--(v25);\draw[e](v15)--(v25);\draw[e](v16)--(v25);\draw[e](v17)--(v26);\draw[e](v18)--(v26);\draw[e](v19)--(v26);\draw[e](v20)--(v26);\draw[e](v21)--(v27);\draw[e](v22)--(v27);\draw[e](v23)--(v27);\draw[e](v24)--(v27);\draw[e](v28)--(v27);\draw[e](v28)--(v26);\draw[e](v25)--(v28);
\end{tikzpicture}
\caption{Iteration 1.1 : obtention of the coloring $\pi_{0}'$}
\label{2}
\end{center}
\end{figure}

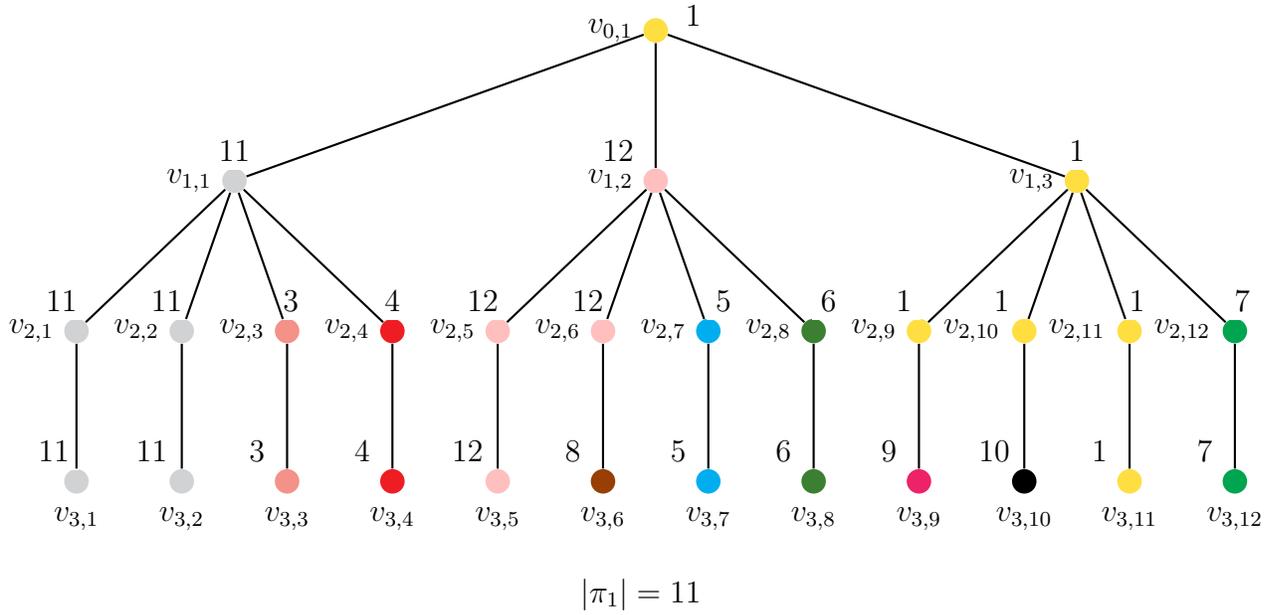
\begin{figure}[ht!]
\begin{center}
\begin{tikzpicture}[inner sep=1.15mm]
\tikzstyle{a}=[circle,fill=Periwinkle!]
\tikzstyle{A}=[circle,fill=Salmon!]
\tikzstyle{L}=[circle,fill=YellowOrange!]
\tikzstyle{b}=[circle,fill=Red!100]
\tikzstyle{o}=[circle,draw,fill=White!]
\tikzstyle{M}=[circle,fill=Dandelion!]
\tikzstyle{O}=[circle,fill=Goldenrod!]
\tikzstyle{S}=[circle,fill=Orange!100]
\tikzstyle{m}=[circle,fill=SeaGreen!]
\tikzstyle{d}=[circle,fill=OliveGreen!]
\tikzstyle{D}=[circle,fill=WildStrawberry!]
\tikzstyle{Q}=[circle,fill=LimeGreen!]
\tikzstyle{N}=[circle,fill=Thistle!]
\tikzstyle{p}=[circle,fill=Fuchsia!]
\tikzstyle{P}=[circle,fill=RoyalPurple!]
\tikzstyle{R}=[circle,fill=Magenta!]
\tikzstyle{I}=[circle,fill=Yellow!100]
\tikzstyle{i}=[circle,fill=GreenYellow!]
\tikzstyle{K}=[circle,fill=CarnationPink!]
\tikzstyle{w}=[circle,fill=Bittersweet!]
\tikzstyle{E}=[circle,fill=black!100]
\tikzstyle{n}=[rectangle,fill=black!0]
\tikzstyle{h}=[circle,fill=blue!100]
\tikzstyle{s}=[circle,fill=TealBlue!]
\tikzstyle{H}=[circle,fill=RoyalBlue!]
\tikzstyle{C}=[circle,fill=Cyan!]
\tikzstyle{f}=[circle,fill=Green!]
\tikzstyle{g}=[circle,fill=SpringGreen!]
\tikzstyle{G}=[circle,fill=Tan!]
\tikzstyle{r}=[circle,fill=Gray!]
\tikzstyle{c}=[circle,fill=CadetBlue!]
\tikzstyle{j}=[circle,fill=pink!]
\tikzstyle{k}=[circle,fill=Peach!]
\tikzstyle{q}=[circle,fill=JungleGreen!]
\tikzstyle{t}=[circle,fill=Apricot!]
\tikzstyle{l}=[circle,fill=BrickRed!]
\tikzstyle{B}=[circle,fill=RawSienna!]
\tikzstyle{J}=[circle,fill=Turquoise!40]
\tikzstyle{F}=[circle,fill=Black!20]
\tikzstyle{T}=[circle,fill=Black!40]
\tikzstyle{u}=[circle,fill=BrickRed!40]
\tikzstyle{U}=[circle,fill=black!50]
\tikzstyle{1}=[circle,fill=brown!100]
\tikzstyle{2}=[circle,fill=purple!100]
\tikzstyle{3}=[circle,fill=green!50]
\tikzstyle{4}=[circle,fill=orange!]
\tikzstyle{e}=[-,thick]

\node [F](v1)at (0,0){};\node [n](v1111)at (0,-0.5){$v_{3,1}$};
\node [F](v2)at (1.4,0){};\node [n](v2222)at (1.4,-0.5){$v_{3,2}$};
\node [A](v3)at (2.8,0){};\node [n](v3333)at (2.8,-0.5){$v_{3,3}$};
\node [b](v4)at (4.2,0){};\node [n](v4444)at (4.2,-0.5){$v_{3,4}$};
\node [j](v5)at (5.6,0){};\node [n](v5555)at (5.6,-0.5){$v_{3,5}$};
\node [B](v6)at (7,0){};\node [n](v6666)at (7,-0.5){$v_{3,6}$};
\node [C](v7)at (8.4,0){};\node [n](v7777)at (8.4,-0.5){$v_{3,7}$};
\node [d](v8)at (9.8,0){};\node [n](v8888)at (9.8,-0.5){$v_{3,8}$};
\node [D](v9)at (11.2,0){};\node [n](v9999)at (11.2,-0.5){$v_{3,9}$};
\node [E](v10)at (12.6,0){};\node [n](v101010)at (12.6,-0.5){$v_{3,10}$};
\node [O](v11)at (14,0){};\node [n](v111111)at (14,-0.5){$v_{3,11}$};
\node [f](v12)at (15.4,0){};\node [n](v121212)at (15.4,-0.5){$v_{3,12}$};
\node [F](v13)at (0,2){};\node [n](v131313)at (-0.6,2){$v_{2,1}$};
\node [F](v14)at (1.4,2){};\node [n](v1414)at (0.8,2){$v_{2,2}$};
\node [A](v15)at (2.8,2){};\node [n](v1515)at (2.2,2){$v_{2,3}$};
\node [b](v16)at (4.2,2){};\node [n](v1616)at (3.6,2){$v_{2,4}$};
\node [j](v17)at (5.6,2){};\node [n](v1717)at (5,2){$v_{2,5}$};
\node [j](v18)at (7,2){};\node [n](v1818)at (6.4,2){$v_{2,6}$};
\node [C](v19)at (8.4,2){};\node [n](v1919)at (7.8,2){$v_{2,7}$};
\node [d](v20)at (9.8,2){};\node [n](v2020)at (9.2,2){$v_{2,8}$};
\node [O](v21)at (11.2,2){};\node [n](v2121)at (10.6,2){$v_{2,9}$};
\node [O](v22)at (12.6,2){};\node [n](v22222)at (11.9,2){$v_{2,10}$};
\node [O](v23)at (14,2){};\node [n](v2323)at (13.3,2){$v_{2,11}$};
\node [f](v24)at (15.4,2){};\node [n](v2424)at (14.7,2){$v_{2,12}$};
\node [F](v25)at (2.1,4){};\node [n](v2525)at (1.5,4){$v_{1,1}$};
\node [j](v26)at (7.7,4){};\node [n](v2626)at (7.1,4){$v_{1,2}$};
\node [O](v27)at (13.3,4){};\node [n](v2727)at (12.7,4){$v_{1,3}$};
\node [O](v28)at (7.7,6){};\node [n](v2828)at (7.1,6){$v_{0,1}$};
\node [n](v41)at (7.5,-1.5){$|\pi_{1}|=11$};
\node [n](v11111)at (-0.3,0.4){$11$};
\node [n](v22222)at (1,0.4){$11$};
\node [n](v33333)at (2.4,0.4){$3$};
\node [n](v44444)at (3.8,0.4){$4$};
\node [n](v55555)at (5.2,0.4){$12$};
\node [n](v66666)at (6.6,0.4){$8$};
\node [n](v77777)at (8,0.4){$5$};
\node [n](v88888)at (9.4,0.4){$6$};
\node [n](v99999)at (10.8,0.4){$9$};
\node [n](v10101010)at (12.2,0.4){$10$};
\node [n](v11111111)at (13.6,0.4){$1$};
\node [n](v12121212)at (15,0.4){$7$};
\node [n](v13131313)at (-0.2,2.4){$11$};
\node [n](v141414)at (1.2,2.4){$11$};
\node [n](v151515)at (2.85,2.4){$3$};
\node [n](v161616)at (4.2,2.4){$4$};
\node [n](v171717)at (5.4,2.4){$12$};
\node [n](v181818)at (6.8,2.4){$12$};
\node [n](v191919)at (8.6,2.4){$5$};
\node [n](v202020)at (10,2.4){$6$};
\node [n](v212121)at (11,2.4){$1$};
\node [n](v2222222)at (12.3,2.4){$1$};
\node [n](v232323)at (14.1,2.4){$1$};
\node [n](v242424)at (15.5,2.4){$7$};
\node [n](v252525)at (2.1,4.4){$11$};
\node [n](v262626)at (7.2,4.4){$12$};
\node [n](v272728)at (13.3,4.4){$1$};
\node [n](v282828)at (8.2,6.2){$1$};
\draw[e](v1)--(v13);\draw[e](v2)--(v14);\draw[e](v3)--(v15);
\draw[e](v4)--(v16);\draw[e](v5)--(v17);\draw[e](v6)--(v18);\draw[e](v7)--(v19);\draw[e](v8)--(v20);\draw[e](v9)--(v21);\draw[e](v10)--(v22);\draw[e](v11)--(v23);\draw[e](v12)--(v24);\draw[e](v13)--(v25);\draw[e](v14)--(v25);\draw[e](v15)--(v25);\draw[e](v16)--(v25);\draw[e](v17)--(v26);\draw[e](v18)--(v26);\draw[e](v19)--(v26);\draw[e](v20)--(v26);\draw[e](v21)--(v27);\draw[e](v22)--(v27);\draw[e](v23)--(v27);\draw[e](v24)--(v27);\draw[e](v28)--(v27);\draw[e](v28)--(v26);\draw[e](v25)--(v28);
\end{tikzpicture}
\caption{Iteration 1.2 : obtention of the quorum coloring $\pi_{1}$}
\label{3}
\end{center}
\end{figure}

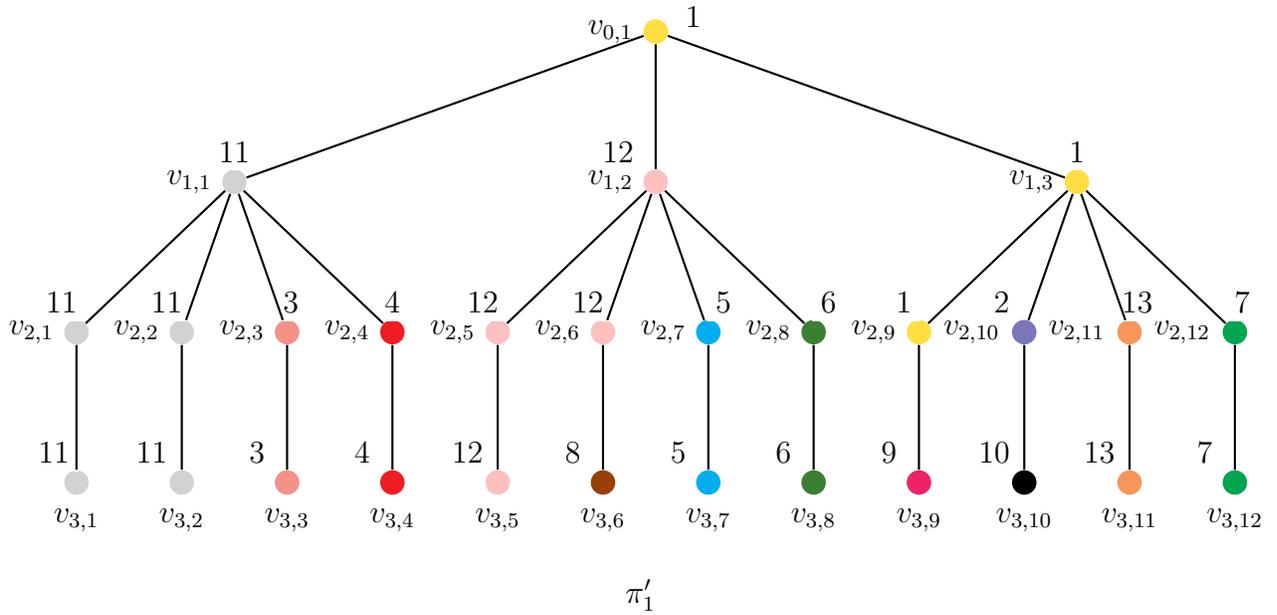
\begin{figure}[ht!]
\begin{center}
\begin{tikzpicture}[inner sep=1.15mm]
\tikzstyle{a}=[circle,fill=Periwinkle!]
\tikzstyle{A}=[circle,fill=Salmon!]
\tikzstyle{L}=[circle,fill=YellowOrange!]
\tikzstyle{b}=[circle,fill=Red!100]
\tikzstyle{o}=[circle,draw,fill=White!]
\tikzstyle{M}=[circle,fill=Dandelion!]
\tikzstyle{O}=[circle,fill=Goldenrod!]
\tikzstyle{S}=[circle,fill=Orange!100]
\tikzstyle{m}=[circle,fill=SeaGreen!]
\tikzstyle{d}=[circle,fill=OliveGreen!]
\tikzstyle{D}=[circle,fill=WildStrawberry!]
\tikzstyle{Q}=[circle,fill=LimeGreen!]
\tikzstyle{N}=[circle,fill=Thistle!]
\tikzstyle{p}=[circle,fill=Fuchsia!]
\tikzstyle{P}=[circle,fill=RoyalPurple!]
\tikzstyle{R}=[circle,fill=Magenta!]
\tikzstyle{I}=[circle,fill=Yellow!100]
\tikzstyle{i}=[circle,fill=GreenYellow!]
\tikzstyle{K}=[circle,fill=CarnationPink!]
\tikzstyle{w}=[circle,fill=Bittersweet!]
\tikzstyle{E}=[circle,fill=black!100]
\tikzstyle{n}=[rectangle,fill=black!0]
\tikzstyle{h}=[circle,fill=blue!100]
\tikzstyle{s}=[circle,fill=TealBlue!]
\tikzstyle{H}=[circle,fill=RoyalBlue!]
\tikzstyle{C}=[circle,fill=Cyan!]
\tikzstyle{f}=[circle,fill=Green!]
\tikzstyle{g}=[circle,fill=SpringGreen!]
\tikzstyle{G}=[circle,fill=Tan!]
\tikzstyle{r}=[circle,fill=Gray!]
\tikzstyle{c}=[circle,fill=CadetBlue!]
\tikzstyle{j}=[circle,fill=pink!]
\tikzstyle{k}=[circle,fill=Peach!]
\tikzstyle{q}=[circle,fill=JungleGreen!]
\tikzstyle{t}=[circle,fill=Apricot!]
\tikzstyle{l}=[circle,fill=BrickRed!]
\tikzstyle{B}=[circle,fill=RawSienna!]
\tikzstyle{J}=[circle,fill=Turquoise!40]
\tikzstyle{F}=[circle,fill=Black!20]
\tikzstyle{T}=[circle,fill=Black!40]
\tikzstyle{u}=[circle,fill=BrickRed!40]
\tikzstyle{U}=[circle,fill=black!50]
\tikzstyle{1}=[circle,fill=brown!100]
\tikzstyle{2}=[circle,fill=purple!100]
\tikzstyle{3}=[circle,fill=green!50]
\tikzstyle{4}=[circle,fill=orange!]
\tikzstyle{e}=[-,thick]
\node [F](v1)at (0,0){};\node [n](v1111)at (0,-0.5){$v_{3,1}$};
\node [F](v2)at (1.4,0){};\node [n](v2222)at (1.4,-0.5){$v_{3,2}$};
\node [A](v3)at (2.8,0){};\node [n](v3333)at (2.8,-0.5){$v_{3,3}$};
\node [b](v4)at (4.2,0){};\node [n](v4444)at (4.2,-0.5){$v_{3,4}$};
\node [j](v5)at (5.6,0){};\node [n](v5555)at (5.6,-0.5){$v_{3,5}$};
\node [B](v6)at (7,0){};\node [n](v6666)at (7,-0.5){$v_{3,6}$};
\node [C](v7)at (8.4,0){};\node [n](v7777)at (8.4,-0.5){$v_{3,7}$};
\node [d](v8)at (9.8,0){};\node [n](v8888)at (9.8,-0.5){$v_{3,8}$};
\node [D](v9)at (11.2,0){};\node [n](v9999)at (11.2,-0.5){$v_{3,9}$};
\node [E](v10)at (12.6,0){};\node [n](v101010)at (12.6,-0.5){$v_{3,10}$};
\node [k](v11)at (14,0){};\node [n](v111111)at (14,-0.5){$v_{3,11}$};
\node [f](v12)at (15.4,0){};\node [n](v121212)at (15.4,-0.5){$v_{3,12}$};
\node [F](v13)at (0,2){};\node [n](v131313)at (-0.6,2){$v_{2,1}$};
\node [F](v14)at (1.4,2){};\node [n](v1414)at (0.8,2){$v_{2,2}$};
\node [A](v15)at (2.8,2){};\node [n](v1515)at (2.2,2){$v_{2,3}$};
\node [b](v16)at (4.2,2){};\node [n](v1616)at (3.6,2){$v_{2,4}$};
\node [j](v17)at (5.6,2){};\node [n](v1717)at (5,2){$v_{2,5}$};
\node [j](v18)at (7,2){};\node [n](v1818)at (6.4,2){$v_{2,6}$};
\node [C](v19)at (8.4,2){};\node [n](v1919)at (7.8,2){$v_{2,7}$};
\node [d](v20)at (9.8,2){};\node [n](v2020)at (9.2,2){$v_{2,8}$};
\node [O](v21)at (11.2,2){};\node [n](v2121)at (10.6,2){$v_{2,9}$};
\node [a](v22)at (12.6,2){};\node [n](v22222)at (11.9,2){$v_{2,10}$};
\node [k](v23)at (14,2){};\node [n](v2323)at (13.3,2){$v_{2,11}$};
\node [f](v24)at (15.4,2){};\node [n](v2424)at (14.7,2){$v_{2,12}$};
\node [F](v25)at (2.1,4){};\node [n](v2525)at (1.5,4){$v_{1,1}$};
\node [j](v26)at (7.7,4){};\node [n](v2626)at (7.1,4){$v_{1,2}$};
\node [O](v27)at (13.3,4){};\node [n](v2727)at (12.7,4){$v_{1,3}$};
\node [O](v28)at (7.7,6){};\node [n](v2828)at (7.1,6){$v_{0,1}$};
\node [n](v41)at (7.5,-1.5){$\pi'_{1}$};
\node [n](v11111)at (-0.3,0.4){$11$};
\node [n](v22222)at (1,0.4){$11$};
\node [n](v33333)at (2.4,0.4){$3$};
\node [n](v44444)at (3.8,0.4){$4$};
\node [n](v55555)at (5.2,0.4){$12$};
\node [n](v66666)at (6.6,0.4){$8$};
\node [n](v77777)at (8,0.4){$5$};
\node [n](v88888)at (9.4,0.4){$6$};
\node [n](v99999)at (10.8,0.4){$9$};
\node [n](v10101010)at (12.2,0.4){$10$};
\node [n](v11111111)at (13.6,0.4){$13$};
\node [n](v12121212)at (15,0.4){$7$};
\node [n](v13131313)at (-0.2,2.4){$11$};
\node [n](v141414)at (1.2,2.4){$11$};
\node [n](v151515)at (2.85,2.4){$3$};
\node [n](v161616)at (4.2,2.4){$4$};
\node [n](v171717)at (5.4,2.4){$12$};
\node [n](v181818)at (6.8,2.4){$12$};
\node [n](v191919)at (8.6,2.4){$5$};
\node [n](v202020)at (10,2.4){$6$};
\node [n](v212121)at (11,2.4){$1$};
\node [n](v2222222)at (12.3,2.4){$2$};
\node [n](v232323)at (14.1,2.4){$13$};
\node [n](v242424)at (15.5,2.4){$7$};
\node [n](v252525)at (2.1,4.4){$11$};
\node [n](v262626)at (7.2,4.4){$12$};
\node [n](v272728)at (13.3,4.4){$1$};
\node [n](v282828)at (8.2,6.2){$1$};
\draw[e](v1)--(v13);\draw[e](v2)--(v14);\draw[e](v3)--(v15);
\draw[e](v4)--(v16);\draw[e](v5)--(v17);\draw[e](v6)--(v18);\draw[e](v7)--(v19);\draw[e](v8)--(v20);\draw[e](v9)--(v21);\draw[e](v10)--(v22);\draw[e](v11)--(v23);\draw[e](v12)--(v24);\draw[e](v13)--(v25);\draw[e](v14)--(v25);\draw[e](v15)--(v25);\draw[e](v16)--(v25);\draw[e](v17)--(v26);\draw[e](v18)--(v26);\draw[e](v19)--(v26);\draw[e](v20)--(v26);\draw[e](v21)--(v27);\draw[e](v22)--(v27);\draw[e](v23)--(v27);\draw[e](v24)--(v27);\draw[e](v28)--(v27);\draw[e](v28)--(v26);\draw[e](v25)--(v28);
\end{tikzpicture}
\caption{Iteration 2.1 : obtention of the coloring $\pi_{1}'$}
\label{4}
\end{center}
\end{figure}

\begin{figure}[ht!]
\begin{center}
\begin{tikzpicture}[inner sep=1.15mm]
\tikzstyle{a}=[circle,fill=Periwinkle!]
\tikzstyle{A}=[circle,fill=Salmon!]
\tikzstyle{L}=[circle,fill=YellowOrange!]
\tikzstyle{b}=[circle,fill=Red!100]
\tikzstyle{o}=[circle,draw,fill=White!]
\tikzstyle{M}=[circle,fill=Dandelion!]
\tikzstyle{O}=[circle,fill=Goldenrod!]
\tikzstyle{S}=[circle,fill=Orange!100]
\tikzstyle{m}=[circle,fill=SeaGreen!]
\tikzstyle{d}=[circle,fill=OliveGreen!]
\tikzstyle{D}=[circle,fill=WildStrawberry!]
\tikzstyle{Q}=[circle,fill=LimeGreen!]
\tikzstyle{N}=[circle,fill=Thistle!]
\tikzstyle{p}=[circle,fill=Fuchsia!]
\tikzstyle{P}=[circle,fill=RoyalPurple!]
\tikzstyle{R}=[circle,fill=Magenta!]
\tikzstyle{I}=[circle,fill=Yellow!100]
\tikzstyle{i}=[circle,fill=GreenYellow!]
\tikzstyle{K}=[circle,fill=CarnationPink!]
\tikzstyle{w}=[circle,fill=Bittersweet!]
\tikzstyle{E}=[circle,fill=black!100]
\tikzstyle{n}=[rectangle,fill=black!0]
\tikzstyle{h}=[circle,fill=blue!100]
\tikzstyle{s}=[circle,fill=TealBlue!]
\tikzstyle{H}=[circle,fill=RoyalBlue!]
\tikzstyle{C}=[circle,fill=Cyan!]
\tikzstyle{f}=[circle,fill=Green!]
\tikzstyle{g}=[circle,fill=SpringGreen!]
\tikzstyle{G}=[circle,fill=Tan!]
\tikzstyle{r}=[circle,fill=Gray!]
\tikzstyle{c}=[circle,fill=CadetBlue!]
\tikzstyle{j}=[circle,fill=pink!]
\tikzstyle{k}=[circle,fill=Peach!]
\tikzstyle{q}=[circle,fill=JungleGreen!]
\tikzstyle{t}=[circle,fill=Apricot!]
\tikzstyle{l}=[circle,fill=BrickRed!]
\tikzstyle{B}=[circle,fill=RawSienna!]
\tikzstyle{J}=[circle,fill=Turquoise!40]
\tikzstyle{F}=[circle,fill=Black!20]
\tikzstyle{T}=[circle,fill=Black!40]
\tikzstyle{u}=[circle,fill=BrickRed!40]
\tikzstyle{U}=[circle,fill=black!50]
\tikzstyle{1}=[circle,fill=brown!100]
\tikzstyle{2}=[circle,fill=purple!100]
\tikzstyle{3}=[circle,fill=green!50]
\tikzstyle{4}=[circle,fill=orange!]
\tikzstyle{e}=[-,thick]
\node [F](v1)at (0,0){};\node [n](v1111)at (0,-0.5){$v_{3,1}$};
\node [F](v2)at (1.4,0){};\node [n](v2222)at (1.4,-0.5){$v_{3,2}$};
\node [A](v3)at (2.8,0){};\node [n](v3333)at (2.8,-0.5){$v_{3,3}$};
\node [b](v4)at (4.2,0){};\node [n](v4444)at (4.2,-0.5){$v_{3,4}$};
\node [j](v5)at (5.6,0){};\node [n](v5555)at (5.6,-0.5){$v_{3,5}$};
\node [B](v6)at (7,0){};\node [n](v6666)at (7,-0.5){$v_{3,6}$};
\node [C](v7)at (8.4,0){};\node [n](v7777)at (8.4,-0.5){$v_{3,7}$};
\node [d](v8)at (9.8,0){};\node [n](v8888)at (9.8,-0.5){$v_{3,8}$};
\node [D](v9)at (11.2,0){};\node [n](v9999)at (11.2,-0.5){$v_{3,9}$};
\node [a](v10)at (12.6,0){};\node [n](v101010)at (12.6,-0.5){$v_{3,10}$};
\node [k](v11)at (14,0){};\node [n](v111111)at (14,-0.5){$v_{3,11}$};
\node [f](v12)at (15.4,0){};\node [n](v121212)at (15.4,-0.5){$v_{3,12}$};
\node [F](v13)at (0,2){};\node [n](v131313)at (-0.6,2){$v_{2,1}$};
\node [F](v14)at (1.4,2){};\node [n](v1414)at (0.8,2){$v_{2,2}$};
\node [A](v15)at (2.8,2){};\node [n](v1515)at (2.2,2){$v_{2,3}$};
\node [b](v16)at (4.2,2){};\node [n](v1616)at (3.6,2){$v_{2,4}$};
\node [j](v17)at (5.6,2){};\node [n](v1717)at (5,2){$v_{2,5}$};
\node [j](v18)at (7,2){};\node [n](v1818)at (6.4,2){$v_{2,6}$};
\node [C](v19)at (8.4,2){};\node [n](v1919)at (7.8,2){$v_{2,7}$};
\node [d](v20)at (9.8,2){};\node [n](v2020)at (9.2,2){$v_{2,8}$};
\node [O](v21)at (11.2,2){};\node [n](v2121)at (10.6,2){$v_{2,9}$};
\node [a](v22)at (12.6,2){};\node [n](v22222)at (11.9,2){$v_{2,10}$};
\node [k](v23)at (14,2){};\node [n](v2323)at (13.3,2){$v_{2,11}$};
\node [f](v24)at (15.4,2){};\node [n](v2424)at (14.7,2){$v_{2,12}$};
\node [F](v25)at (2.1,4){};\node [n](v2525)at (1.5,4){$v_{1,1}$};
\node [j](v26)at (7.7,4){};\node [n](v2626)at (7.1,4){$v_{1,2}$};
\node [O](v27)at (13.3,4){};\node [n](v2727)at (12.7,4){$v_{1,3}$};
\node [O](v28)at (7.7,6){};\node [n](v2828)at (7.1,6){$v_{0,1}$};
\node [n](v41)at (7.5,-1.5){$|\pi_{2}|=12$};
\node [n](v11111)at (-0.3,0.4){$11$};
\node [n](v22222)at (1,0.4){$11$};
\node [n](v33333)at (2.4,0.4){$3$};
\node [n](v44444)at (3.8,0.4){$4$};
\node [n](v55555)at (5.2,0.4){$12$};
\node [n](v66666)at (6.6,0.4){$8$};
\node [n](v77777)at (8,0.4){$5$};
\node [n](v88888)at (9.4,0.4){$6$};
\node [n](v99999)at (10.8,0.4){$9$};
\node [n](v10101010)at (12.2,0.4){$2$};
\node [n](v11111111)at (13.6,0.4){$13$};
\node [n](v12121212)at (15,0.4){$7$};
\node [n](v13131313)at (-0.2,2.4){$11$};
\node [n](v141414)at (1.2,2.4){$11$};
\node [n](v151515)at (2.85,2.4){$3$};
\node [n](v161616)at (4.2,2.4){$4$};
\node [n](v171717)at (5.4,2.4){$12$};
\node [n](v181818)at (6.8,2.4){$12$};
\node [n](v191919)at (8.6,2.4){$5$};
\node [n](v202020)at (10,2.4){$6$};
\node [n](v212121)at (11,2.4){$1$};
\node [n](v2222222)at (12.3,2.4){$2$};
\node [n](v232323)at (14.1,2.4){$13$};
\node [n](v242424)at (15.5,2.4){$7$};
\node [n](v252525)at (2.1,4.4){$11$};
\node [n](v262626)at (7.2,4.4){$12$};
\node [n](v272728)at (13.3,4.4){$1$};
\node [n](v282828)at (8.2,6.2){$1$};
\draw[e](v1)--(v13);\draw[e](v2)--(v14);\draw[e](v3)--(v15);
\draw[e](v4)--(v16);\draw[e](v5)--(v17);\draw[e](v6)--(v18);\draw[e](v7)--(v19);\draw[e](v8)--(v20);\draw[e](v9)--(v21);\draw[e](v10)--(v22);\draw[e](v11)--(v23);\draw[e](v12)--(v24);\draw[e](v13)--(v25);\draw[e](v14)--(v25);\draw[e](v15)--(v25);\draw[e](v16)--(v25);\draw[e](v17)--(v26);\draw[e](v18)--(v26);\draw[e](v19)--(v26);\draw[e](v20)--(v26);\draw[e](v21)--(v27);\draw[e](v22)--(v27);\draw[e](v23)--(v27);\draw[e](v24)--(v27);\draw[e](v28)--(v27);\draw[e](v28)--(v26);\draw[e](v25)--(v28);
\end{tikzpicture}
\caption{Iteration 2.2 : obtention of the quorum coloring $\pi_{2}$}
\label{5}
\end{center}
\end{figure}
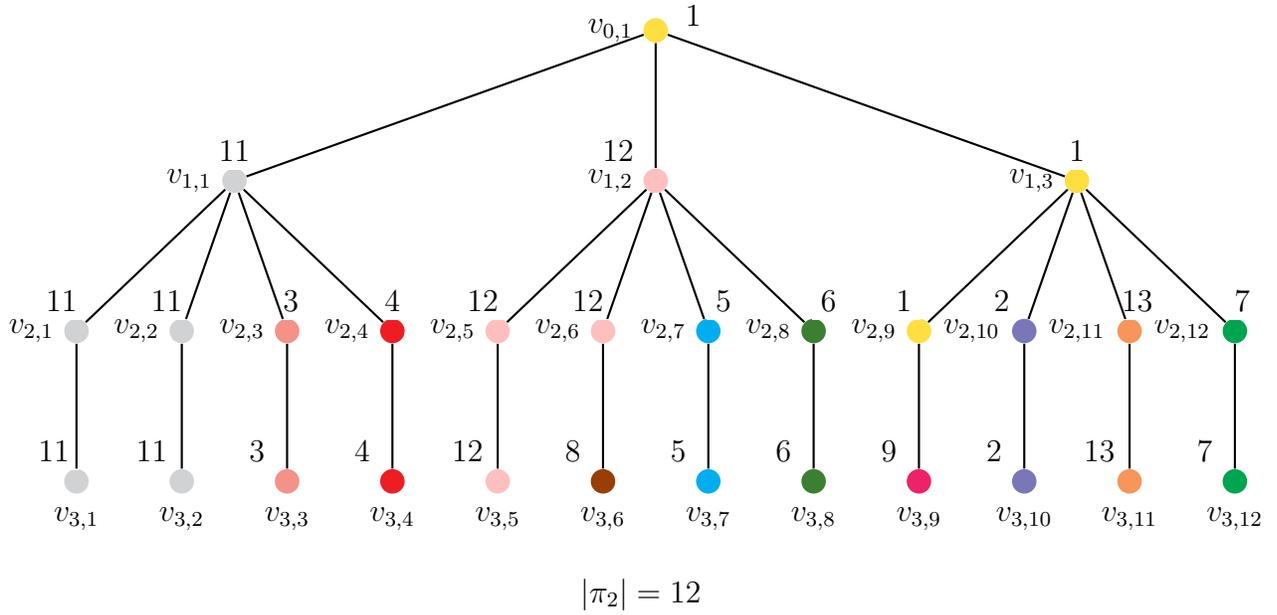

\begin{figure}[ht!]
\begin{center}
\begin{tikzpicture}[inner sep=1.15mm]
\tikzstyle{a}=[circle,fill=Periwinkle!]
\tikzstyle{A}=[circle,fill=Salmon!]
\tikzstyle{L}=[circle,fill=YellowOrange!]
\tikzstyle{b}=[circle,fill=Red!100]
\tikzstyle{o}=[circle,draw,fill=White!]
\tikzstyle{M}=[circle,fill=Dandelion!]
\tikzstyle{O}=[circle,fill=Goldenrod!]
\tikzstyle{S}=[circle,fill=Orange!100]
\tikzstyle{m}=[circle,fill=SeaGreen!]
\tikzstyle{d}=[circle,fill=OliveGreen!]
\tikzstyle{D}=[circle,fill=WildStrawberry!]
\tikzstyle{Q}=[circle,fill=LimeGreen!]
\tikzstyle{N}=[circle,fill=Thistle!]
\tikzstyle{p}=[circle,fill=Fuchsia!]
\tikzstyle{P}=[circle,fill=RoyalPurple!]
\tikzstyle{R}=[circle,fill=Magenta!]
\tikzstyle{I}=[circle,fill=Yellow!100]
\tikzstyle{i}=[circle,fill=GreenYellow!]
\tikzstyle{K}=[circle,fill=CarnationPink!]
\tikzstyle{w}=[circle,fill=Bittersweet!]
\tikzstyle{E}=[circle,fill=black!100]
\tikzstyle{n}=[rectangle,fill=black!0]
\tikzstyle{h}=[circle,fill=blue!100]
\tikzstyle{s}=[circle,fill=TealBlue!]
\tikzstyle{H}=[circle,fill=RoyalBlue!]
\tikzstyle{C}=[circle,fill=Cyan!]
\tikzstyle{f}=[circle,fill=Green!]
\tikzstyle{g}=[circle,fill=SpringGreen!]
\tikzstyle{G}=[circle,fill=Tan!]
\tikzstyle{r}=[circle,fill=Gray!]
\tikzstyle{c}=[circle,fill=CadetBlue!]
\tikzstyle{j}=[circle,fill=pink!]
\tikzstyle{k}=[circle,fill=Peach!]
\tikzstyle{q}=[circle,fill=JungleGreen!]
\tikzstyle{t}=[circle,fill=Apricot!]
\tikzstyle{l}=[circle,fill=BrickRed!]
\tikzstyle{B}=[circle,fill=RawSienna!]
\tikzstyle{J}=[circle,fill=Turquoise!40]
\tikzstyle{F}=[circle,fill=Black!20]
\tikzstyle{T}=[circle,fill=Black!40]
\tikzstyle{u}=[circle,fill=BrickRed!40]
\tikzstyle{U}=[circle,fill=black!50]
\tikzstyle{1}=[circle,fill=brown!100]
\tikzstyle{2}=[circle,fill=purple!100]
\tikzstyle{3}=[circle,fill=green!50]
\tikzstyle{4}=[circle,fill=orange!]
\tikzstyle{e}=[-,thick]
\node [E](v1)at (0,0){};\node [n](v1111)at (0,-0.5){$v_{3,1}$};
\node [K](v2)at (1.4,0){};\node [n](v2222)at (1.4,-0.5){$v_{3,2}$};
\node [A](v3)at (2.8,0){};\node [n](v3333)at (2.8,-0.5){$v_{3,3}$};
\node [b](v4)at (4.2,0){};\node [n](v4444)at (4.2,-0.5){$v_{3,4}$};
\node [p](v5)at (5.6,0){};\node [n](v5555)at (5.6,-0.5){$v_{3,5}$};
\node [B](v6)at (7,0){};\node [n](v6666)at (7,-0.5){$v_{3,6}$};
\node [C](v7)at (8.4,0){};\node [n](v7777)at (8.4,-0.5){$v_{3,7}$};
\node [d](v8)at (9.8,0){};\node [n](v8888)at (9.8,-0.5){$v_{3,8}$};
\node [D](v9)at (11.2,0){};\node [n](v9999)at (11.2,-0.5){$v_{3,9}$};
\node [a](v10)at (12.6,0){};\node [n](v101010)at (12.6,-0.5){$v_{3,10}$};
\node [k](v11)at (14,0){};\node [n](v111111)at (14,-0.5){$v_{3,11}$};
\node [f](v12)at (15.4,0){};\node [n](v121212)at (15.4,-0.5){$v_{3,12}$};
\node [F](v13)at (0,2){};\node [n](v131313)at (-0.6,2){$v_{2,1}$};
\node [F](v14)at (1.4,2){};\node [n](v1414)at (0.8,2){$v_{2,2}$};
\node [A](v15)at (2.8,2){};\node [n](v1515)at (2.2,2){$v_{2,3}$};
\node [b](v16)at (4.2,2){};\node [n](v1616)at (3.6,2){$v_{2,4}$};
\node [j](v17)at (5.6,2){};\node [n](v1717)at (5,2){$v_{2,5}$};
\node [j](v18)at (7,2){};\node [n](v1818)at (6.4,2){$v_{2,6}$};
\node [C](v19)at (8.4,2){};\node [n](v1919)at (7.8,2){$v_{2,7}$};
\node [d](v20)at (9.8,2){};\node [n](v2020)at (9.2,2){$v_{2,8}$};
\node [O](v21)at (11.2,2){};\node [n](v2121)at (10.6,2){$v_{2,9}$};
\node [a](v22)at (12.6,2){};\node [n](v22222)at (11.9,2){$v_{2,10}$};
\node [k](v23)at (14,2){};\node [n](v2323)at (13.3,2){$v_{2,11}$};
\node [f](v24)at (15.4,2){};\node [n](v2424)at (14.7,2){$v_{2,12}$};
\node [F](v25)at (2.1,4){};\node [n](v2525)at (1.5,4){$v_{1,1}$};
\node [j](v26)at (7.7,4){};\node [n](v2626)at (7.1,4){$v_{1,2}$};
\node [O](v27)at (13.3,4){};\node [n](v2727)at (12.7,4){$v_{1,3}$};
\node [O](v28)at (7.7,6){};\node [n](v2828)at (7.1,6){$v_{0,1}$};
\node [n](v41)at (7.5,-1.5){$\pi'_{2}$};
\node [n](v11111)at (-0.3,0.4){$10$};
\node [n](v22222)at (1,0.4){$14$};
\node [n](v33333)at (2.4,0.4){$3$};
\node [n](v44444)at (3.8,0.4){$4$};
\node [n](v55555)at (5.2,0.4){$15$};
\node [n](v66666)at (6.6,0.4){$8$};
\node [n](v77777)at (8,0.4){$5$};
\node [n](v88888)at (9.4,0.4){$6$};
\node [n](v99999)at (10.8,0.4){$9$};
\node [n](v10101010)at (12.2,0.4){$2$};
\node [n](v11111111)at (13.6,0.4){$13$};
\node [n](v12121212)at (15,0.4){$7$};
\node [n](v13131313)at (-0.2,2.4){$11$};
\node [n](v141414)at (1.2,2.4){$11$};
\node [n](v151515)at (2.85,2.4){$3$};
\node [n](v161616)at (4.2,2.4){$4$};
\node [n](v171717)at (5.4,2.4){$12$};
\node [n](v181818)at (6.8,2.4){$12$};
\node [n](v191919)at (8.6,2.4){$5$};
\node [n](v202020)at (10,2.4){$6$};
\node [n](v212121)at (11,2.4){$1$};
\node [n](v2222222)at (12.3,2.4){$2$};
\node [n](v232323)at (14.1,2.4){$13$};
\node [n](v242424)at (15.5,2.4){$7$};
\node [n](v252525)at (2.1,4.4){$11$};
\node [n](v262626)at (7.2,4.4){$12$};
\node [n](v272728)at (13.3,4.4){$1$};
\node [n](v282828)at (8.2,6.2){$1$};
\draw[e](v1)--(v13);\draw[e](v2)--(v14);\draw[e](v3)--(v15);
\draw[e](v4)--(v16);\draw[e](v5)--(v17);\draw[e](v6)--(v18);\draw[e](v7)--(v19);\draw[e](v8)--(v20);\draw[e](v9)--(v21);\draw[e](v10)--(v22);\draw[e](v11)--(v23);\draw[e](v12)--(v24);\draw[e](v13)--(v25);\draw[e](v14)--(v25);\draw[e](v15)--(v25);\draw[e](v16)--(v25);\draw[e](v17)--(v26);\draw[e](v18)--(v26);\draw[e](v19)--(v26);\draw[e](v20)--(v26);\draw[e](v21)--(v27);\draw[e](v22)--(v27);\draw[e](v23)--(v27);\draw[e](v24)--(v27);\draw[e](v28)--(v27);\draw[e](v28)--(v26);\draw[e](v25)--(v28);
\end{tikzpicture}
\caption{Iteration 3.1 : obtention of the coloring $\pi_{2}'$}
\label{6}
\end{center}
\end{figure}
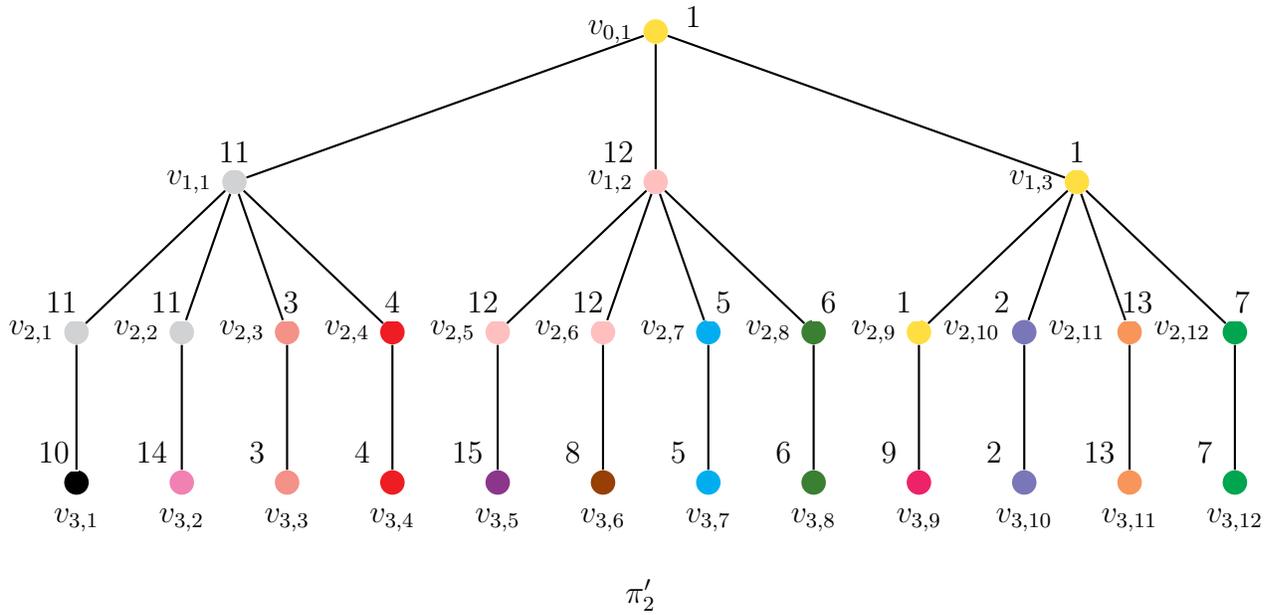

\begin{figure}[ht!]
\begin{center}
\begin{tikzpicture}[inner sep=1.15mm]
\tikzstyle{a}=[circle,fill=Periwinkle!]
\tikzstyle{A}=[circle,fill=Salmon!]
\tikzstyle{L}=[circle,fill=YellowOrange!]
\tikzstyle{b}=[circle,fill=Red!100]
\tikzstyle{o}=[circle,draw,fill=White!]
\tikzstyle{M}=[circle,fill=Dandelion!]
\tikzstyle{O}=[circle,fill=Goldenrod!]
\tikzstyle{S}=[circle,fill=Orange!100]
\tikzstyle{m}=[circle,fill=SeaGreen!]
\tikzstyle{d}=[circle,fill=OliveGreen!]
\tikzstyle{D}=[circle,fill=WildStrawberry!]
\tikzstyle{Q}=[circle,fill=LimeGreen!]
\tikzstyle{N}=[circle,fill=Thistle!]
\tikzstyle{p}=[circle,fill=Fuchsia!]
\tikzstyle{P}=[circle,fill=RoyalPurple!]
\tikzstyle{R}=[circle,fill=Magenta!]
\tikzstyle{I}=[circle,fill=Yellow!100]
\tikzstyle{i}=[circle,fill=GreenYellow!]
\tikzstyle{K}=[circle,fill=CarnationPink!]
\tikzstyle{w}=[circle,fill=Bittersweet!]
\tikzstyle{E}=[circle,fill=black!100]
\tikzstyle{n}=[rectangle,fill=black!0]
\tikzstyle{h}=[circle,fill=blue!100]
\tikzstyle{s}=[circle,fill=TealBlue!]
\tikzstyle{H}=[circle,fill=RoyalBlue!]
\tikzstyle{C}=[circle,fill=Cyan!]
\tikzstyle{f}=[circle,fill=Green!]
\tikzstyle{g}=[circle,fill=SpringGreen!]
\tikzstyle{G}=[circle,fill=Tan!]
\tikzstyle{r}=[circle,fill=Gray!]
\tikzstyle{c}=[circle,fill=CadetBlue!]
\tikzstyle{j}=[circle,fill=pink!]
\tikzstyle{k}=[circle,fill=Peach!]
\tikzstyle{q}=[circle,fill=JungleGreen!]
\tikzstyle{t}=[circle,fill=Apricot!]
\tikzstyle{l}=[circle,fill=BrickRed!]
\tikzstyle{B}=[circle,fill=RawSienna!]
\tikzstyle{J}=[circle,fill=Turquoise!40]
\tikzstyle{F}=[circle,fill=Black!20]
\tikzstyle{T}=[circle,fill=Black!40]
\tikzstyle{u}=[circle,fill=BrickRed!40]
\tikzstyle{U}=[circle,fill=black!50]
\tikzstyle{1}=[circle,fill=brown!100]
\tikzstyle{2}=[circle,fill=purple!100]
\tikzstyle{3}=[circle,fill=green!50]
\tikzstyle{4}=[circle,fill=orange!]
\tikzstyle{e}=[-,thick]
\node [E](v1)at (0,0){};\node [n](v1111)at (0,-0.5){$v_{3,1}$};
\node [K](v2)at (1.4,0){};\node [n](v2222)at (1.4,-0.5){$v_{3,2}$};
\node [A](v3)at (2.8,0){};\node [n](v3333)at (2.8,-0.5){$v_{3,3}$};
\node [b](v4)at (4.2,0){};\node [n](v4444)at (4.2,-0.5){$v_{3,4}$};
\node [p](v5)at (5.6,0){};\node [n](v5555)at (5.6,-0.5){$v_{3,5}$};
\node [B](v6)at (7,0){};\node [n](v6666)at (7,-0.5){$v_{3,6}$};
\node [C](v7)at (8.4,0){};\node [n](v7777)at (8.4,-0.5){$v_{3,7}$};
\node [d](v8)at (9.8,0){};\node [n](v8888)at (9.8,-0.5){$v_{3,8}$};
\node [D](v9)at (11.2,0){};\node [n](v9999)at (11.2,-0.5){$v_{3,9}$};
\node [a](v10)at (12.6,0){};\node [n](v101010)at (12.6,-0.5){$v_{3,10}$};
\node [k](v11)at (14,0){};\node [n](v111111)at (14,-0.5){$v_{3,11}$};
\node [f](v12)at (15.4,0){};\node [n](v121212)at (15.4,-0.5){$v_{3,12}$};
\node [F](v13)at (0,2){};\node [n](v131313)at (-0.6,2){$v_{2,1}$};
\node [F](v14)at (1.4,2){};\node [n](v1414)at (0.8,2){$v_{2,2}$};
\node [A](v15)at (2.8,2){};\node [n](v1515)at (2.2,2){$v_{2,3}$};
\node [b](v16)at (4.2,2){};\node [n](v1616)at (3.6,2){$v_{2,4}$};
\node [j](v17)at (5.6,2){};\node [n](v1717)at (5,2){$v_{2,5}$};
\node [j](v18)at (7,2){};\node [n](v1818)at (6.4,2){$v_{2,6}$};
\node [C](v19)at (8.4,2){};\node [n](v1919)at (7.8,2){$v_{2,7}$};
\node [d](v20)at (9.8,2){};\node [n](v2020)at (9.2,2){$v_{2,8}$};
\node [O](v21)at (11.2,2){};\node [n](v2121)at (10.6,2){$v_{2,9}$};
\node [a](v22)at (12.6,2){};\node [n](v22222)at (11.9,2){$v_{2,10}$};
\node [k](v23)at (14,2){};\node [n](v2323)at (13.3,2){$v_{2,11}$};
\node [f](v24)at (15.4,2){};\node [n](v2424)at (14.7,2){$v_{2,12}$};
\node [F](v25)at (2.1,4){};\node [n](v2525)at (1.5,4){$v_{1,1}$};
\node [j](v26)at (7.7,4){};\node [n](v2626)at (7.1,4){$v_{1,2}$};
\node [O](v27)at (13.3,4){};\node [n](v2727)at (12.7,4){$v_{1,3}$};
\node [O](v28)at (7.7,6){};\node [n](v2828)at (7.1,6){$v_{0,1}$};
\node [n](v11111)at (-0.3,0.4){$10$};
\node [n](v22222)at (1,0.4){$14$};
\node [n](v33333)at (2.4,0.4){$3$};
\node [n](v44444)at (3.8,0.4){$4$};
\node [n](v55555)at (5.2,0.4){$15$};
\node [n](v66666)at (6.6,0.4){$8$};
\node [n](v77777)at (8,0.4){$5$};
\node [n](v88888)at (9.4,0.4){$6$};
\node [n](v99999)at (10.8,0.4){$9$};
\node [n](v10101010)at (12.2,0.4){$2$};
\node [n](v11111111)at (13.6,0.4){$13$};
\node [n](v12121212)at (15,0.4){$7$};
\node [n](v13131313)at (-0.2,2.4){$11$};
\node [n](v141414)at (1.2,2.4){$11$};
\node [n](v151515)at (2.85,2.4){$3$};
\node [n](v161616)at (4.2,2.4){$4$};
\node [n](v171717)at (5.4,2.4){$12$};
\node [n](v181818)at (6.8,2.4){$12$};
\node [n](v191919)at (8.6,2.4){$5$};
\node [n](v202020)at (10,2.4){$6$};
\node [n](v212121)at (11,2.4){$1$};
\node [n](v2222222)at (12.3,2.4){$2$};
\node [n](v232323)at (14.1,2.4){$13$};
\node [n](v242424)at (15.5,2.4){$7$};
\node [n](v252525)at (2.1,4.4){$11$};
\node [n](v262626)at (7.2,4.4){$12$};
\node [n](v272728)at (13.3,4.4){$1$};
\node [n](v282828)at (8.2,6.2){$1$};
\node [n](v41)at (7.5,-1.5){$|\pi_{3}|=15>|\pi_{0}|=10$};
\draw[e](v1)--(v13);\draw[e](v2)--(v14);\draw[e](v3)--(v15);
\draw[e](v4)--(v16);\draw[e](v5)--(v17);\draw[e](v6)--(v18);\draw[e](v7)--(v19);\draw[e](v8)--(v20);\draw[e](v9)--(v21);\draw[e](v10)--(v22);\draw[e](v11)--(v23);\draw[e](v12)--(v24);\draw[e](v13)--(v25);\draw[e](v14)--(v25);\draw[e](v15)--(v25);\draw[e](v16)--(v25);\draw[e](v17)--(v26);\draw[e](v18)--(v26);\draw[e](v19)--(v26);\draw[e](v20)--(v26);\draw[e](v21)--(v27);\draw[e](v22)--(v27);\draw[e](v23)--(v27);\draw[e](v24)--(v27);\draw[e](v28)--(v27);\draw[e](v28)--(v26);\draw[e](v25)--(v28);
\end{tikzpicture}
\caption{Iteration 3.2 : obtention of the quorum coloring $\pi_{3}$}
\label{7}
\end{center}
\end{figure}
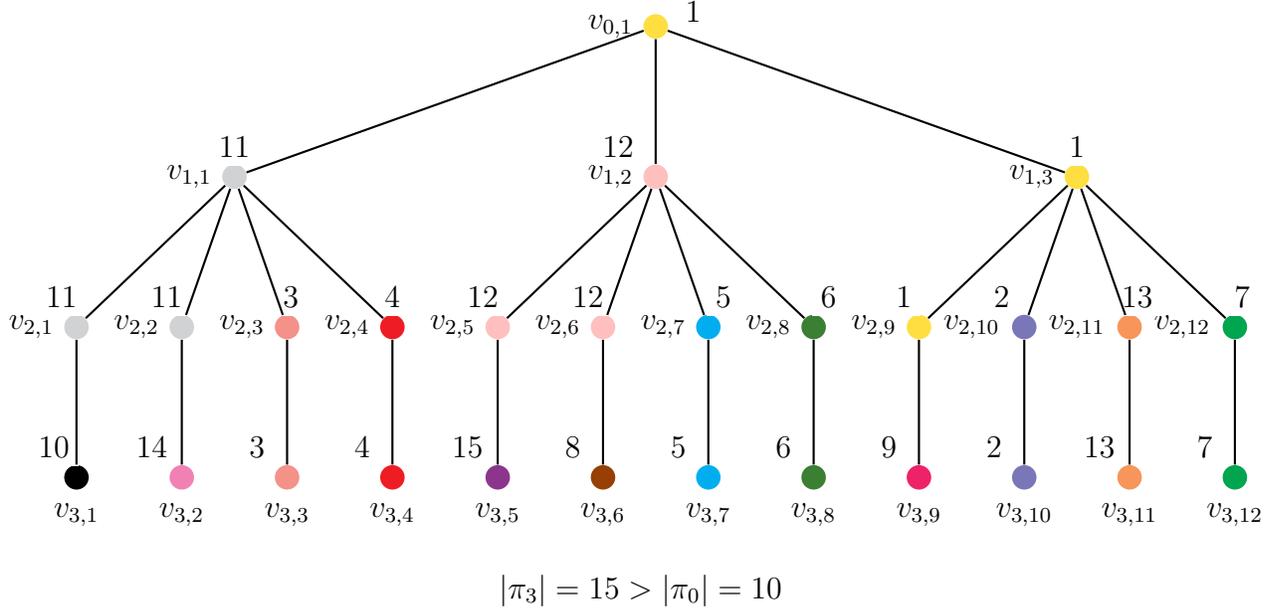

Figures~\ref{1}-\ref{7} illustrate how Algorithm 1 works on a perfect $N_{i}$-ary tree per level with $\displaystyle\left(N_{i}\right)_{0\leq i\leq 2}=(3,4,1)$.

\vspace{0.2cm}

\vspace{0.2cm}

\newpage

\vspace{0.75cm}

\vspace{0.75cm}

\newpage

\vspace{0.75cm}

\vspace{0.75cm}

\newpage

\vspace{0.4cm}

We are now ready to prove our central result by designing a linear-time algorithm both for finding a $\psi_{q}$-coloring and computing the quorum coloring number of any perfect $N_{i}$-ary tree per level.

\begin{theorem}\label{linear}
For every perfect $N_{i}$-ary tree per level $T$ of order $n$, one can both find a cost-effective $\psi_{q}$-coloring and compute $\psi_{q}(T)$ in $O(n)$.
\end{theorem}

\begin{proof} Let $T=(V,E,r)$ be a perfect $N_{i}$-ary tree per level of order $n$. We run the following algorithm denoted by Algorithm 2.

\begin{itemize} 

\item[1.] Assign a color to $r$ and go to step 2.

\item[2.] Assign arbitrarily the $r$'s color to exactly $\left\lfloor\dfrac{d_{T}(r)}{2}\right\rfloor$ of its children, color its $\left\lceil\dfrac{d_{T}(r)}{2}\right\rceil$ remaining children with $\left\lceil\dfrac{d_{T}(r)}{2}\right\rceil$ new colors and set $\alpha_{1}^{1}=\left\lceil\dfrac{N_{0}}{2}\right\rceil+1$. Then, go to step 3.

\item[3.] For $i=1$ to $h-1$ and $j=1$ to $\ell_{i}$, consider the following mutually exclusive two cases. 

\begin{description}

\item[(i)] If $v_{i,j}$ has the same color as its parent, then assign arbitrarily the color of $v_{i,j}$ to exactly $\left\lceil\dfrac{N_{i}}{2}\right\rceil-1$ of its children, color its $\left\lfloor\dfrac{N_{i}}{2}\right\rfloor+1$ remaining others with $\left\lfloor\dfrac{N_{i}}{2}\right\rfloor+1$ new colors and set $\alpha_{i+1}^{j}=\left\lfloor\dfrac{N_{i}}{2}\right\rfloor+1$.

\item[(ii)] Otherwise, if $v_{i,j}$ has not the same color as its parent, then assign arbitrarily the $v_{i,j}$'s color to exactly $\left\lceil\dfrac{N_{i}}{2}\right\rceil$ of its children, color its $\left\lfloor\dfrac{N_{i}}{2}\right\rfloor$ remaining children with $\left\lfloor\dfrac{N_{i}}{2}\right\rfloor$ new colors and set $\alpha_{i+1}^{j}=\left\lfloor\dfrac{N_{i}}{2}\right\rfloor$.

\end{description}

\item[4.] Set $\displaystyle\alpha=\sum_{i=0}^{h-1}\sum_{j=1}^{\ell_{i}}\alpha_{i+1}^{j}$.

\end{itemize}

Algorithm 2 must terminate since $T$ is finite. By Assertion 2 of Observation \ref{Obs6}, this algorithm concludes by producing a cost-effective quorum coloring $\pi$ of $T$ as output thanks to steps 2 and 3. Moreover, since the $D_{i}$'s vertices have the same degree and therefore play the same role in $T$ (indeed, their descendants that are at the same distance all have the same degree), then by removing the vertex labeling "`$v_{i,j}$"' at the end of the run of Algorithm 2, one can easily see that 
 there exists a unique cost-effective quorum coloring of $T$ up to isomorphism. 
 It follows by Corollary \ref{corQC->CEF} that $\pi$ is 
 a $\psi_{q}$- coloring of $T$. 
  Furthermore, on the one hand the total number of comparison tests between the color of a vertex and that of its parent is upper bounded by the order $n$ of $T$. On the other hand, the total number of color assignments is upper bounded by $n$ too. In addition, the computation of the $\alpha_{i+1}^{j}$'s requires at most $2n$ operations. Finally, the sum calculated in step 4 is performed in at most $n$ operations. Consequently, the total number of elementary operations necessary to both find $\pi$ and calculate $\psi_{q}(T)=\alpha$ is at most equal to $n+n+2n+n< 5n$, hence the theorem. \end{proof}

The next four corollaries follow immediately from Theorem \ref{linear} and from the fact that every perfect $N$-ary tree is also a perfect $N_{i}$-ary tree per level with $N_{i}=N$ for every $i\in\{0,1,\ldots,h-1\}$.

\begin{corollary}\label{Cor1Linear} For every perfect $N_{i}$-ary tree per level $T$ of order $n$, one can find a cost-effective $\psi_{q}$-coloring in $O(n)$.\end{corollary}

\begin{corollary}\label{Cor2Linear} For every perfect $N_{i}$-ary tree per level $T$ of order $n$, one can compute $\psi_{q}(T)$ in $O(n)$.\end{corollary}

\begin{corollary}\label{Cor3Linear} For every perfect $N$-ary tree $T$ of order $n$, one can find a cost-effective $\psi_{q}$-coloring in $O(n)$.\end{corollary}

\begin{corollary}\label{Cor4Linear} For every perfect $N$-ary tree $T$ of order $n$, one can compute $\psi_{q}(T)$ in $O(n)$.\end{corollary}

It can be seen without difficulty that the quorum coloring $\pi_{3}$ of Figure~\ref{7} can be obtained by Algorithm 2, which illustrates the unicity of a $\psi_{q}$-coloring in a perfect $N_{i}$-ary tree per level. In particular, one can check that we have : 

\[\alpha_{1}^{1}=\left\lceil\dfrac{N_{0}}{2}\right\rceil+1=\left\lceil\dfrac{3}{2}\right\rceil+1=3,~\alpha_{2}^{1}=\left\lfloor\dfrac{N_{1}}{2}\right\rfloor = \left\lfloor\dfrac{4}{2}\right\rfloor=2,~\alpha_{2}^{2}=\left\lfloor\dfrac{N_{1}}{2}\right\rfloor = \left\lfloor\dfrac{4}{2}\right\rfloor=2,\]

\[\alpha_{2}^{3}=\left\lfloor\dfrac{N_{1}}{2}\right\rfloor+1=\left\lfloor\dfrac{4}{2}\right\rfloor+1=3,~\alpha_{3}^{1}=\left\lfloor\dfrac{N_{2}}{2}\right\rfloor+1=\left\lfloor\dfrac{1}{2}\right\rfloor+1=1,~\alpha_{3}^{2}=\left\lfloor\dfrac{N_{2}}{2}\right\rfloor+1=\left\lfloor\dfrac{1}{2}\right\rfloor+1=1,\]

\[\alpha_{3}^{3}=\left\lfloor\dfrac{N_{2}}{2}\right\rfloor = \left\lfloor\dfrac{1}{2}\right\rfloor=0,~\alpha_{3}^{4}=\left\lfloor\dfrac{N_{2}}{2}\right\rfloor = \left\lfloor\dfrac{1}{2}\right\rfloor=0,~\alpha_{3}^{7}=\left\lfloor\dfrac{N_{2}}{2}\right\rfloor = \left\lfloor\dfrac{1}{2}\right\rfloor=0,\]

\[\alpha_{3}^{8}=\left\lfloor\dfrac{N_{2}}{2}\right\rfloor = \left\lfloor\dfrac{1}{2}\right\rfloor=0,~\alpha_{3}^{5}=\left\lfloor\dfrac{N_{2}}{2}\right\rfloor+1=\left\lfloor\dfrac{1}{2}\right\rfloor+1=1,~\alpha_{3}^{6}=\left\lfloor\dfrac{N_{2}}{2}\right\rfloor+1=\left\lfloor\dfrac{1}{2}\right\rfloor+1=1,\]

\[\alpha_{3}^{9}=\left\lfloor\dfrac{N_{2}}{2}\right\rfloor+1=\left\lfloor\dfrac{1}{2}\right\rfloor+1=1,~\alpha_{3}^{10}=\left\lfloor\dfrac{N_{2}}{2}\right\rfloor = \left\lfloor\dfrac{1}{2}\right\rfloor=0,~\alpha_{3}^{11}=\left\lfloor\dfrac{N_{2}}{2}\right\rfloor = \left\lfloor\dfrac{1}{2}\right\rfloor=0,\]

\[\alpha_{3}^{12}=\left\lfloor\dfrac{N_{2}}{2}\right\rfloor = \left\lfloor\dfrac{1}{2}\right\rfloor=0\text{ and that }\displaystyle\alpha=\sum_{i=0}^{h-1}\sum_{j=1}^{\ell_{i}}\alpha_{i+1}^{j}=15.\]

\section{Concluding remarks}

In this article, we have designed a linear-time algorithm finding a quorum coloring of maximum cardinality together with the quorum coloring number of any perfect tree satisfying the property that the vertices at the same depth have the same degree. Since every perfect $N$-ary tree is a perfect $N_{i}$-ary tree per level, then the same algorithm is applicable for perfect $N$-ary trees, which solves the second part of Problem 2. This being established, there remain three questions of interest to go further in our investigation. First, one can think of generalizing somewhat the work done in this paper to perfect trees whose siblings have the same degree; we will call them {\it locally perfect $N_{i,j}$-ary trees}. The children of a given vertex would then not necessarily play the same role since the descendants of any two siblings of same depth would not necessarily have the same degree. 
 Which leads us to our first question.







\begin{description}

\item[{\bf 3.}] Is Algorithm 2 adaptable to locally perfect $N_{i,j}$-ary trees ? If yes, how would we choose the siblings taking the color of their parent (or equivalently of those that do not take it) in step 3?
    
\end{description}

More generally, we can study the following question.

\begin{description}

\item[{\bf 4.}] Can you design a linear-time algorithm for finding $\psi_{q}(T)$ for any perfect trees $T$?
    
\end{description}

The study of Questions 3 and 4 could constitute two steps leading to the resolution of Question 1.

Finally, the third and last question is combinatorial and is stated as follows.

\begin{description}

\item[{\bf 5.}] Is it possible to determine the exact value of the quorum coloring number of a perfect $N$-ary tree or a perfect $N_{i}$-ary tree by level ?

\end{description}

For a possible inductive approach, one can verify without difficulty using the exact value of Corollary \ref{Cor0} that for any perfect binary tree $T$ of order $n$ and height $h$, we have as initialization of the induction that \[\displaystyle\psi_{q}(T)=\dfrac{2^{h+2}}{2^{ 2}-1} \left[1-\left(\dfrac{1}{4}\right)^{ \left\lfloor\frac{h}{2}\right\rfloor+1}\right].\]
















\end{document}